\documentclass[notitlepage, 11pt]{article}
\usepackage{graphicx}

\usepackage{bbm}

\usepackage{color,xcolor}
\usepackage{latexsym, soul}
\usepackage{cancel}
\usepackage{amsmath, amssymb, amsfonts}
\usepackage[toc]{appendix}
\usepackage[normalem]{ulem}

\numberwithin{equation}{section}

\DeclareMathOperator*{\argmax}{arg\,max}

\newtheorem{lemma}{Lemma}
\newtheorem{theorem}{Theorem}

\newtheorem{proposition}{Proposition}

\newtheorem{assumption}{Assumption}

\newtheorem{remark}{Remark}

\newcommand{\beginsec}{
\setcounter{lemma}{0}
\setcounter{theorem}{0}
\setcounter{corollary}{0}
\setcounter{definition}{0}
\setcounter{example}{0}
\setcounter{proposition}{0}
\setcounter{condition}{0}
\setcounter{assumption}{0}
\setcounter{conjecture}{0}
\setcounter{problem}{0}
\setcounter{remark}{0}
}

\newcommand{\noi}{\noindent}

\newcommand{\E}{\mathbb{E}}
\newcommand{\R}{\mathbb{R}}

\newcommand{\N}{\mathbb{N}}

\newcommand{\la}{\lambda}

\newcommand{\eps}{\varepsilon}

\newcommand{\one}{\mathbbm{1}}

\newcommand{\Q}{{\mathbb Q}}

\newcommand{\PP}{{\mathbb P}}

\newcommand{\brho}{\boldsymbol{\rho}}

\newcommand{\calA}{{\cal A}}

\newcommand{\calC}{{\cal C}}
\newcommand{\calD}{{\cal D}}
\newcommand{\calE}{{\cal E}}
\newcommand{\calF}{{\cal F}}
\newcommand{\calG}{{\cal G}}

\newcommand{\calL}{{\cal L}}

\newcommand{\calO}{{\cal O}}

\newcommand{\calQ}{{\cal Q}}

\newcommand{\calX}{{\cal X}}

\newcommand{\skp}{\vspace{\baselineskip}}
\newcommand{\supp}{{\rm supp}}

\newcommand{\To}{\Rightarrow}

\newcommand\iy{\infty}

\newcommand{\qed}{\hfill $\Box$}

\newcommand{\limn}{\lim_{n\to\iy}}

\newcommand{\clc}{{\cal C}}

\newcommand{\clo}{{\cal O}}

\newcommand{\cld}{{\cal D}}

\newcommand\redsout{\bgroup\markoverwith{\textcolor{red}{\rule[0.5ex]{2pt}{0.4pt}}}\ULon}

\newcommand{\sredm}[1]{\ifmmode\text{\xout{\ensuremath{\displaystyle \textcolor{red}{#1}}}}\else\sout{\textcolor{red}{#1}}\fi}

\title{Asymptotic optimality of the generalized $c\mu$ rule under model uncertainty\thanks{This is the final version of the paper. To appear in {\it Stochastic Processes and their Applications}.}}
\author{
Asaf Cohen\thanks{Department of Mathematics,
University of Michigan,
Ann Arbor, MI 48109, USA. Email: shloshim@gmail.com, web: https://sites.google.com/site/asafcohentau/ The research of A. Cohen is supported by the National Science Foundation (DMS-2006305).
}\,\,\,
and\,\,\,
Subhamay Saha\thanks{Department of Mathematics,
Indian Institute of Technology Guwahati,
Guwahati, Assam 781039, India. Email: saha.subhamay@iitg.ac.in, web: https://sites.google.com/site/sahasubhamay86/}
}

\date{\today}

\begin{document}

\maketitle

\begin{abstract}
We consider a critically-loaded multiclass queueing control problem with model uncertainty. The model consists of $I$ types of customers and a single server. At any time instant, a decision-maker (DM) allocates the server's effort to the customers. The DM's goal is to minimize a convex holding cost that accounts for the ambiguity with respect to the model, i.e., the arrival and service rates. For this, we consider an adversary player whose role is to choose the worst-case scenario. Specifically, we assume that the DM has a reference probability model in mind and that the cost function is formulated by the supremum over equivalent admissible probability measures to the reference measure with two components, the first is the expected holding cost, and the second one is a penalty for the adversary player for deviating from the reference model. The penalty term is formulated by a general divergence measure. 

We show that although that under the equivalent admissible measures the critically-load condition might be violated, the generalized $c\mu$ rule is asymptotically optimal for this problem.


\skp
\noi{\bf AMS subject classifications.} 60K25, 60J60, 93E20, 60F05, 68M20, 91A15.

\skp
\noi{\bf Keywords:} diffusion scaling, heavy-traffic, model uncertainty, general divergence, generalized c$\mu$ rule.
\end{abstract}

\section{Introduction}\label{sec1}
\beginsec
In this article we consider a multiclass queueing control problem (QCP) under diffusion-scaled heavy-traffic condition, wherein the decision-maker (DM) is uncertain about the underlying model, in the sense that she is unsure about the arrival and service rates. Upon arrival customers are kept in queues in accordance to their types. The DM allocates customers to the server, trying to minimize a convex holding cost functional. The ambiguity modeling is done as follows. We assume that the DM has a reference probability model in mind, and to account for the uncertainty she tries to optimize among a family of models, penalizing the deviation of a model from the reference model using a class of general divergence measures. Thus the cost function the DM is facing is given by
\begin{align}\notag
\sup_{\Q}\Big\{\E^{\Q}\Big[\int_0^\iy \varrho(t)C(\hat  X(t))dt \Big]
 -\calL(\Q\|\PP)\Big\}.
 \end{align}
The supremum is taken over a set of candidate admissible probability measures. The first term is the discounted holding cost term ($\varrho$ is the discounting, $X(t)$ is a vector whose components are the sizes of the queues, and $C$ is a holding cost function). The second term is the penalty term accounting for the ambiguity, where $\PP$ is the reference probability measure, $\Q$ is a candidate probability measure, and $\calL$ is a divergence measure. The complete description of the terms is given in the next section.

In general, QCPs are almost intractable. So, an approach pioneered by Harrison \cite{har1988}, is to approximate the QCP via a diffusion control problem (DCP). To this end, DCP is first solved and then its optimal control is used to find an asymptotically optimal control of the original QCP.\footnote{
Note that here, the pre-limit model (with uncertainty) is of interest. The limiting model is more tractable and from which we are able to derive a simple index policy. Hence, it serves us as an auxiliary model. This is typical, but not restricted to queueing theory; see also the following examples \cite{kim1964} in biology, \cite{igl1969} in actuary, and \cite{coh2015} in parameter estimation. This is in contrast to cases where one may be interested in the limiting model and would use a sequence of approximating models, which are either more tractable or for which numerical analysis is more feasible. When the limiting process is of interest, one may care about stability results, that is: how properties, such as optimal behavior, are affected by (small) misspecifications of the underlying model; see e.g., \cite{lar-zit2007, pap-pos-sap2019, bac-bar-bei2020, bay-dol-guo2020} and the references therein. 
 } In order to carry out the asymptotic analysis, the QCP needs to be  scaled. 
Specifically, 
one considers a sequence of queueing systems, labeled by a scaling parameter $n\in\N$. It is assumed that, for each queue, the {\it traffic intensity}, which is the ratio between the rate of arrivals to the rate of service is $1-\clo(1/\sqrt{n})$ (see \cite{chen-yao}). The arrival and service rates for each queue are of order $\clo(n)$. Queue lengths are scaled (divided) by $\sqrt{n}$. In the literature of queueing systems, this scaling is referred to as the \emph{heavy traffic} diffusion scaling regime, see \cite{bell-will-1, BG2006, BG2012} and references therein. The terminology stems from the fact that given a sequence of `reasonable' controls, one for each $n$, there is a limiting diffusion process, along a subsequence as $n\to\iy$. Note that if the traffic intensity is not of order $1-\clo(1/\sqrt{n})$, then, as $n\to\iy$, the scaled queue length becomes degenerate (converges to either infinity or zero). This explains the terminology \emph{critically loaded} system. 
In absence of the penalizing term (or, alternatively, the supremum), by the martingale central limit theorem argument (see \cite[Section 7]{eth-kur}), the QCP can be approximated by an optimal control problem of a diffusion process. We detail more in the sequel about the specific structure of the QCP studied here, with and without ambiguity.

Aside from the diffusion scaling, there are other well-studied asymptotic regimes for QCPs. The first example we provide is the {\it fluid scaling}, which emerges from the functional law of large numbers and leads to a limiting deterministic control problem, see \cite{Kaspi15, Shimkin14, Walsh14} and references therein. As with the relationship between the central limit theorem and the law of large numbers, the fluid scaling is less fine. Specifically, the traffic intensity is $\clo(1)$ and queues are scaled by $n$. Two other regimes are the the heavy traffic moderate deviation scaling (see \cite{ata-bis, ata-coh, ata-coh2017, bis2014} and references therein) and the large deviation scaling (see \cite{AGS1, AGS2} and references therein).
The optimization criteria in the last two setups are risk-sensitive costs. While the standard diffusion scaling (without ambiguity) leads to a stochastic control problem, the moderate-deviation and large-deviation regimes lead to two-player minimax games, where the minimizer stands for the DM from the QCP and the maximizer is an adverse player that models a worst-case scenario. These games are deterministic. In this paper we follow a similar philosophy and use a diffusion scaling. Since the optimization is over a class of models, the limiting control problem in our case is in fact a stochastic game. 

In QCPs one of the aims is to come up with easily implementable asymptotically optimal policies. This issue is even more crucial when the dimension of the problem is high. One of the most classical such policies in the setup of multiclass  queueing network is the $c\mu$ rule which is asymptotically optimal in the case of a linear holding cost. A generalized version of this was introduced by van Mieghem \cite{van}, the dynamic priority rule known as the generalized $c\mu$ rule where the parameter $c$ is variable and obtained by feedback from the system's state.  Specifically, if the holding cost rate of a class-$i$ customer is given by $C_i$, where $C_i$ are given smooth, convex functions then the rule is to prioritize the classes according to the index
$\mu_iC'_i(\hat Q_i^n(t))$, where $n$ is the scaling parameter,
$\mu_i$ and $\hat Q_i^n(t)$ are the corresponding service rate and
{\it diffusion scaled} queue length at time $t$, and $C'_i$ denotes the derivative of $C_i$. Hence, it is clear that this index policy handles well the curse of dimensionality.
In \cite{ata-sah-cmu}, Atar and Saha show the asymptotic optimality of the generalized $c\mu$ rule in the moderate deviation heavy-traffic regime. However, in the large deviation regime, asymptotic optimality of the generalized $c\mu$ rule is not to be expected. In fact, in \cite{AGS2}, the authors show that a rule other than the $c\mu$ rule is optimal in the setup of linear holding cost.

Thus the main contribution of this article is that we extend the robustness of the generalized $c\mu$ rule as an asymptotically optimal policy to the setup of model uncertainty with a general divergence measure. Model uncertainty is a very realistic assumption in the real life situation. Due to the complexity of real-world systems, lack of sufficient calibration, and inaccurate assumptions,
one cannot precisely model the arrival and departure processes, see e.g. \cite{Hasenbein17, ward15, whitt2006staffing} and references therein. In recent years there has been increasing interest in robust analysis of stochastic control problems. We consider uncertainty in the diffusion scale of the QCP in a way that is also referred to as the Knightian uncertainty with a reference model.\footnote{Another type of Knightian uncertainty that is considered in the literature is one where there is no reference model and the worst case scenario is formulated via a collection of measures, which are not necessarily equivalent to each other. This type of ambiguity gained popularity in the mathematical finance community and is associated with the G-expectation introduced by Peng \cite{pen2007}, for further reading see e.g., \cite{pos-gui-tou, neu-nut, ban-dol-gok} and the references therein.} 
As described in the first paragraph of the introduction, the DM chooses a stochastic control and the maximizer responses in a non-anticipating way with respect to the proper filtration by choosing a `worst case scenario'. A popular choice for the penalty function $\calL(\Q\|\PP)$ is the {\it Kullback--Leibler divergence}, which due to the tractability it yields, has been widely exploited in the literature in recent years in various fields, e.g., in 
control theory \cite{pet-jam-dup2000, iye2005, nil-el2005, blanchet2014robust},
actuary and finance \cite{maenhout2004robust, gla-xu2014, bay-zha}, 
economics \cite{han-ser2006}, 
distributional robust optimization \cite{ben-den-de-mel-ren2013}, 
and in queueing systems (not under the heavy traffic regime though) \cite{jai-lim-sha2010, blanchet2014robust,  lam2016}. This choice leads to an equivalent stochastic differential game, where the maximizer responses to the DM's control by stochastically perturbing the drift of the state dynamics in a non-anticipating way using a stochastic process and is quadratically penalized.  Uncertainty in fluid models of queueing was studied e.g., in \cite{bassamboo2010capacity, dup2003, whitt2006staffing}. Recently, Krishnasamy et al.~\cite{2018arXiv180206723K} took a different approach from the Knightian uncertainty considered here and studied a learning-based variant of the $c\mu$ rule for scheduling in multi-class queueing systems, where the service rates are unknown and the DM's objective is to minimize a regret criterion comparing between the $c\mu$ rule that uses the empirically learned service rates and the $c\mu$ rule with the known rates.

This paper continues a line of research initiated in \cite{ASAFCOHEN2018, coh2019a}, where asymptotic analysis of a multiclass QCP under uncertainty is being studied with linear holding cost, finite buffers, and a penalty ambiguity term that is formulated via the Kullback--Leibler divergence. Paper \cite{ASAFCOHEN2018} studies the limiting stochastic game and characterizes the value function as the unique solution to a free-boundary ordinary differential equation. The solvability of this differential equation heavily relies on the special structure of the Kullback--Leibler divergence, which is shown to lead to an equivalent quadratic penalty term. Paper \cite{coh2019a} establishes asymptotic optimality using the limiting problem and crucially relies on the ordinary differential equation and again on the quadratic penalty term that follows by the Kullback--Leibler divergence. To compare with the current work, we consider here a convex holding cost, infinite buffers, and a general divergence measure. Thus, the model considered in this paper is very different from \cite{ASAFCOHEN2018, coh2019a}. The more general structure of the model adds to the subtlety of the arguments as we detail in the next paragraph. Another very important contribution of this paper that it gives an insight about a sufficient condition on the structure of the penalty term depending on the growth rate of the cost function, so that asymptotic optimality of the generalized $c\mu$ rule is still true. 

After setting up the limiting stochastic minimax game, the asymptotic optimality of the generalized $c\mu$ rule is established by showing that asymptotically the value function of this game forms both lower and upper bound for the QCP. The important fact that we take advantage of is that the limiting game problem has an explicit solution (in terms of optimal control for the minimizer and value function) in terms of the Skorohod map and a minimizing curve describing the optimal workload distribution in the limiting problem. Therefore, we avoid differential equations analysis. In order to establish the lower bound, we take an arbitrary sequence of controls in the QCP and compare it with the cost in the limiting game associated with the optimal control of the minimizer. Recall that the cost function in both the QCP and the limiting game involve suprema over probability measures. Here, one should pay attention that  the randomness in the limiting problem is generated by a Brownian motion and in the QCP by Poisson random measures. Hence, for the comparison between the suprema we use a discretization technique due to Kushner, see Lemma \ref{lem41}. For this, we associate each admissible probability measure in the limiting problem with a stochastic process via the Radon--Nykodim derivative. Then, we approximate it by stochastic processes adapted to the filtrations of the QCPs, which in turn are translated to admissible probability measures in the QCP; these approximated processes are the arrival and service rates under the admissible measures. In addition, the proof of the lower bound uses tightness and martingale arguments under the approximated probability measure.  
For the proof of the upper bound we equip the minimizer of the QCP with the generalized $c\mu$ rule and fix an arbitrary sequence of admissible probability measures for the maximizer.  Again, the probability measures are translated to rates. We show in two steps that the maximizer will abstain from making large perturbations to the arrival and service rates, in the sense that in average the critically load heavy-traffic condition is preserved, see Proposition \ref{prop42}. Then, we use a truncation technique to approximate the rates with bounded ones, see Proposition \ref{prop43}. The rest of the proof uses a state-collapse property and again tightness and martingale arguments. 

The rest of the paper is organized as follows. In the next subsection we enlist the notations used throughout the paper. In Section 2, we describe the queueing model and the robust optimization problem. In Section 3, we describe the limiting stochastic game and its solution. In section 4, we describe the generalized $c\mu$ rule and prove its asymptotic optimality by proving the convergence of the prelimit value functions to the limiting value function. 
\subsection{Notation}
We use the following notation.
For $a,b\in\R$, $a\wedge b=\min\{a,b\}$ and $a\vee b=\max\{a,b\}$. For a positive integer $k$ and $c,d\in\R^k$, $c\cdot d$ denotes the usual
scalar product and $\|c\|=(c\cdot c)^{1/2}$. We denote $[0, \iy)$ by $\R_+$. For subintervals $I_1,I_2\subseteq\R$ and $m\in\{1,2\}$ we denote by $\clc(I_1,I_2)$, $\clc^m(I_1,I_2)$, and $\cld(I_1, I_2)$ the space of continuous functions [resp., functions with continuous derivatives of order $m$, functions that are right-continuous with finite left limits (RCLL)] mapping $I_1\to I_2$. 
The space $\calD(I_1,I_2)$ is endowed with the usual Skorohod topology.
For $f:\mathbb{R}_+\rightarrow \mathbb{R}^I$ and $t>0$, set $\|f\|_t=\sup_{0\le s\le t}\|f(s)\|$.

\section{The queueing model}\label{sec2}
\beginsec
\subsection{The reference probability space and some fundamental processes}\label{sec21}
\beginsec
Consider a single server model with $I$ classes of customers. Each class has its own designated unbounded buffer. Upon arrival, customers are queued in the corresponding buffers.
Processor sharing is allowed, but two customers from the same class cannot be served
simultaneously. The system is studied under heavy-traffic. Hence, we consider a sequence of systems, indexed by $n\in\N$, which is referred to as the {\it scaling parameter}. Each class of customers is associated with two rates: arrival and service. To account for the ambiguity in each of the total $2I$ parameters it is convenient to work with a probability space in a product form. For every $i\in[I]:=\{1,\ldots,I\}$ and $n\in\N$ we consider two {\it reference probability spaces}
$(\Omega^n_{A,i},\calF^n_{A,i}
\PP^n_{A,i})$ and $(\Omega^n_{S,i},\calF^n_{S,i},
\PP^n_{S,i})$. The first one supports a Poisson process $A^n_i$ with a given rate $\lambda^n_i$ and the second one a Poisson process $S^n_i$ with rate $\mu^n_i$. $A^n_i$ counts the number of arrivals to the $i$-th buffer. The process $S^n_i$ is referred to as the \emph{potential service time} process, in the sense that, for every $t\in\R_+$, $S^n_i(t)$ is the number of service completions of class $i$ customers had the server worked on class $i$ for $t$ units of time.

In order to account for different level of ambiguities for each of the $2I$ processes, we construct the {\it complete reference probability space} that supports the processes $A^n=(A^n_i:i\in[I])$ and $S^n=(S^n_i:i\in[I])$ and which is given in a product form as follows,
\begin{align}\notag
(\Omega^n,\calF^n,
\PP^n):=
\Big(\prod_{i=1}^I(\Omega^n_{A,i}\times\Omega^n_{S,i}),\otimes_{i=1}^I(\calF^n_{A,i}\otimes\calG^n_{S,i})
,\prod_{i=1}^I(\PP^n_{A,i}\times\PP^n_{S,i})\Big),
\end{align}
where $\otimes_{i=1}^I(\calF^n_{A,i}\otimes\calG^n_{S,i})=(\calF^n_{A,1}\otimes\calG^n_{A,2})\otimes\ldots\otimes(\calF^n_{A,I}\otimes\calG^n_{S,I})$.
Notice that from the structure of the probability space it follows that for every fixed $n\in\N$, under the measure $\PP^n$, the processes $A^n_1, S^n_1, \ldots, A^n_I, S^n_I$ are mutually independent. Moreover,
$\PP^n\circ(A^n_i)^{-1}=\PP^n_{A,i}\circ(A^n_i)^{-1}$ and $\PP^n\circ(S^n_i)^{-1}=\PP^n_{S,i}\circ(S^n_i)^{-1}$, $i\in[I]$.

Let $U^n=(U^n_i:i\in[I])$ be an RCLL process taking values in $\{x=(x_1,\ldots,x_I)\in[0,1]^I:\sum x_i\le1\}$. The term $U^n_i(t)$ represents
the fraction of effort devoted at time $t$ by the server to the class-$i$ customer at the head of the line. For each $i\in[I]$, the process  $(T^n_i(t))_{t\in\R_+}$ given by
\begin{align}\label{206}
T^n_i(t):=\int_0^t U^n_i(s)ds,\qquad t\in\R_+,
\end{align}
represents the units of time that the server devoted to class $i$ until time $t$.
For every $t\in\R_+$ and $i\in[I]$, $S^n_i(T^n_i(t))$ is the number of service completions of class $i$ customers until time $t$. This is a Cox process with intensity $\mu^n_iU^n_i$.

Denote by $X^n_i(t)$ the number of class $i$ customers in the system at time $t$. Then,
\begin{equation}\label{204a}
X^n_i(t)=X^n_i(0)+A^n_i(t)-S^n_i(T^n_i(t)),\qquad t\in\R_+,i\in[I].
\end{equation}
The system is assumed to be \emph{critically loaded}. That is, the rate parameters satisfy
\begin{align}\label{205}
\lambda^n_i:=\lambda_i n+\hat\lambda_in^{1/2}+o(n^{1/2}),\qquad
\mu^n_i:=\mu_in+\hat\mu_in^{1/2}+o(n^{1/2}),
\end{align}
where $\lambda_i,\mu_i\in(0,\iy)$ and $\hat\lambda_i,\hat\mu_i\in\R$ are fixed and
$\sum_{i=1}^I\rho_i=1$, where $\rho_i:=\lambda_i/\mu_i$, $i\in[I]$. 
Using the diffusion scaling
\begin{align}\label{207}
\begin{split}
 \hat{A}^n_i(t)&:=n^{-1/2}(A^n_i(t)-\la^n_i t),\quad  \hat{S}^n_i(t):=n^{-1/2}(S^n_i(t)-\mu^n_i t)
,\\
 \hat X^n_i(t)&:=n^{-1/2}X^n_i(t), \quad\hat Y^n_i(t):=\mu^n_in^{-1/2}(\rho_i t-T^n_i(t)),\quad
\hat m^n_i:=n^{-1/2}(\lambda^n_i-\rho_i\mu^n_i),
\end{split}
\end{align}
we obtain the following scaled version of \eqref{204a},
\begin{equation} \notag
 \hat{X}^n_i(t)=\hat{X}^n_i(0) + \hat m^n_it
 +\hat{A}^n_i(t)-\hat{S}^n_i(T^n_i(t))
  +\hat  Y^n_i(t),\qquad t\in\R_+
\end{equation}
Denote $L^n(t)=(L^n_i(t):i\in[I])$ for $L^n=\hat A^n, \hat S^n, T^n, \hat X^n, \hat Y^n$ and also $\hat m^n=(\hat m^n_i:i\in[I])$.
It will be assumed throughout that,
\begin{equation}
  \notag
\exists \limn \hat X^n(0)=: \hat x_0.
\end{equation}
For simplicity, we assume that $\{X^n_i(0)\}_{i,n}$ are deterministic, hence so is $\hat x_0$.

The process $U^n$ is regarded as an {\it admissible control} in the $n$-th system if it is adapted to the filtration ${\calF}^n_t={\calF}^n(t):=\sigma\{A^n_i(s),S^n_i(T^n_i(s)), i\in[I], s\le t\}$ and $U^n_i(t)=0$ whenever $X^n_i(t)=0$.
The latter condition asserts that the server cannot devote any effort to an empty class. We denote the set of admissible controls for the DM in the $n$-th system by $\calA^n$.

\subsection{The robust optimization problem}\label{sec22}
We now describe the cost function.
Fix a discount function $\varrho:\R_+\to\R_+$ and a holding cost $C:\R^I_+\to\R_+$, which satisfy some regularity and growth conditions (see Assumption \ref{assumption1} below).

The {\it risk-neutral optimization problem} is given by
\begin{align}\notag
\inf_{ U^n\in\calA^n}\;\E^{\PP^n}\Big[\int_0^\iy \varrho(t)C(\hat X^n(t))dt\Big].
\end{align}
A variation of this problem with convex delay costs was studied by van Mieghem \cite{van}. 

In order to capture the uncertainty of the DM about the underlying probability measure we consider a set of candidate probability measures and the DM optimizes against the worst measure.
Each such measure is being penalized in accordance to its deviation from the reference measure
by a
divergence.
%
Notice that given a process $\hat Y^n$ adapted to $\calF^n_t$, satisfying \eqref{207}, there exists an admissible control $U^n$ for which \eqref{207} holds with $T^n$ given by \eqref{206}. Hence, we  refer to $\hat Y^n$ as well as the control in the $n$-th system. The DM is facing the following robust optimization problem:
\begin{align}\notag
 V^n:= \;&
\inf_{\hat Y^n\in\calA^n}\;\sup_{\Q^n\in\calQ^n}  J^n(\hat Y^n,\Q^n),
\end{align}
where
\begin{align}\label{212}
 & J^n(\hat Y^n,\Q^n):=\E^{\Q^n}\Big[\int_0^\iy \varrho(t)C(\hat  X^n(t))dt \Big]
 -\sum_{i=1}^I\calL_{A,i}(\Q^n_{A,i}\|\PP^n_{A,i})
 -\sum_{i=1}^I\calL_{S,i}(\Q^n_{S,i}\|\PP^n_{S,i})
 \end{align}
 and its components are the following:
 \begin{itemize}
 \item $\calQ^n$ is
the set of all the measures of the form $\Q^n=\prod_{i=1}^I(\Q^n_{A,i}\times\Q^n_{S,i})$, satisfying
\begin{align}\label{216}
	\frac{d\Q^n_{A,i}}{d\PP^n_{A,i}}(t)&=\exp\Big\{\int_0^t\log\left(\frac{\psi^n_{A,i}(s)}{\la^n_i}\right) d A^n_i(s)-\int_0^t (\psi^n_{A,i}(s)-\la^n_i)ds\Big\}, \\\label{216z}
	\frac{d\Q^n_{S,i}}{d\PP^n_{S,i}}(t)&=\exp\Big\{\int_0^t\log\left(\frac{\psi^n_{S,i}(s)}{\mu^n_i}\right) d S^n_i(T^n_i(s))-\int_0^t (\psi^n_{S,i}(s)-\mu^n_i)dT^n_i(s)\Big\},
\end{align}
for measurable and positive processes $\psi^n_{j,i}$ that are predictable w.r.t.~the filtration generated by the arrival and service completions processes, satisfying $\int_0^t\psi^n_{j,i}(s)ds<\iy$ $\PP^n_{j,i}$-a.s., $j\in\{A,S\}, i\in[I]$. We refer to the elements of $\calQ$ as {\it admissible controls} for the {\it adverse player}, which is also called the {\it maximizer}. Occasionally we  abuse terminology and refer to the processes $\psi^n_{j,i}$ as the maximizer's controls.

\item For every $j\in\{A,S\}$, $i\in[I]$, and equivalent measures $\Q$ and $\PP$, 
the divergence $\calL_{j,i}$, 
is given by
\begin{align}\notag
\calL_{j,i}(\Q\|\PP)&:=\E^{\Q}\Big[\int_0^\iy\varrho(t)g_{j,i}\Big(\log\Big(\frac{d\Q}{d\PP}(t)\Big)\Big)dt\Big]
,
 \end{align}
where $g_{j,i}:\R\to\R$ satisfy some regularity and growth conditions given below in Assumption \ref{assumption1}.
\end{itemize}
\begin{remark}
\begin{enumerate}
\item The conditions asserted on $\{\psi^n_{j,i}\}$ guarantee that the right-hand sides in \eqref{216} are $\PP^n_{j,i}$-martingales, and that under the measure $\Q^n_{A,i}$ (resp., $\Q^n_{S,i}$), the processes $A^n_i$ (resp., $S^n_i(T^n_i)$) is a counting process with infinitesimal intensity $\psi^n_{A,i}(t)dt$ (resp., $\psi^n_{2,i}dT^n_i(t)$). 
\item Notice that we do not assume that the critically load condition is preserved under the measure $\Q^n$, that is $\psi^n_{A,i}(t)-\lambda^n_i=\calO(n^{1/2})$, uniformly over $t\in\R_+$. However, as we show in Proposition \ref{prop42} below, this condition holds in average.
\end{enumerate} 
\end{remark}

%

\begin{assumption}\label{assumption1}
There exist constants $\bar p\ge p\ge 1$, $c_0,c_2>0$, and $c_1\in\R$, such that the functions $C$, $\{g_{j,i}\}_{j\in\{A,S\},i\in[I]}$, and $\varrho$, satisfy the following conditions.
\begin{enumerate}
\item
There are strictly increasing, strictly convex, continuously differentiable functions $C_i:\R_+\to\R_+$, $i\in[I]$, such that,
\begin{itemize}
\item $C_i'(0)=C_i(0)=0$;
\item $C(x)=\sum_{i=1}^IC_i(x_i)$, where $x=(x_i:i\in[I])$;
\item there exist constants $c_0>0$ and $p\ge1$ such that 
for every $i\in[I]$ and $x\in\R_+$, $C_i(x)\le c_0(1+x^p)$.
\end{itemize}
\item For every $j\in\{A,S\}$ and $i\in[I]$, $g_{j,i}$ are convex and non-decreasing. There exist constants $c_1 \in \mathbb{R}$, $c_2,c_3>0$, such that $|g_{j,i}(x)|\leq c_3(1+|x|^{\bar{p}}) $ for any $x\in\R$, and $g_{i,j}(x)\ge c_1+c_2x^{\bar p}  $ for $x\in\R_+$. 
\item $\varrho:\R_+\to\R_+$ is non-increasing and satisfies $\int_0^\iy\varrho(t)t^{\bar p}dt<\iy$.
\end{enumerate}
\end{assumption}
Throughout the entire paper we assume that Assumption \ref{assumption1} is in force.

\begin{remark}\label{rem22}
\begin{enumerate}
\item 
The assumptions on $C_i$, other than the growth condition, are identical to those in \cite{ata-sah-cmu}. They ensure the existence of a continuous minimizing curve which is necessary for the explicit solution of the limiting differential game. The intuition behind the relation between the growth condition on $C_i$ and the lower bound on the growth of $g_{i,j}$ is that it ensures that the maximizer does not perturb the rate by ``too much" and the heavy traffic condition is preserved on the ``average". This suggests the correct order of the penalty term that should be considered if the cost function has more than linear growth. 

\item Note that aside of having finite buffers, the model in 
\cite{ASAFCOHEN2018, coh2019a} considers a linear cost and the KL divergence, which corresponds to the case $p=\bar p=1$. Nevertheless, the linear cost  cannot be considered as a private of the current model because the cost function is strictly convex and its derivative is zero at the origin. A better comparison in this case is the  studied in \cite{ata-sah-cmu} under the moderate-deviation model with risk-sensitive cost. There, the cost term within the exponent satisfies all the assumptions mentioned above and has a linear growth. The penalty term, which emerges from the moderate-deviation scale, corresponds to the case $p=\bar p=1$; that is, $g$ is linear, and the penalty term is expressed in terms of the KL divergence. 

\item The structure of the divergence considered here by composing $g$ on $\log(d\Q/d\PP)$ emerges from the need to set the penalty term and the cost function on the same scale. A similar structure was considered by Hern\'andez-Hern\'andez and Schied in \cite{her-sch} in a control theoretic setup. With an adaptation to our model, the function $g$ is composed on the integrands of \eqref{216}--\eqref{216z}. The cost function in \cite{her-sch} is the log of the state process, which has an exponential growth. Thus, it corresponds to a linear growth of the cost function. In turn, their penalty $g$ has at least a quadratic growth. This is the same growth that the KL divergence has, which corresponds to the linear growth in the cost.

\end{enumerate}
\end{remark}

\section{The limiting problem - a stochastic differential game}\label{sec3}
\beginsec

\subsection{The game setup}\label{sec31}


The QCP is approximated by a two-player stochastic differential game. To set it up we need the following notation.
Set the vectors
\begin{align}\notag
\theta&:=\limn(n/\mu^n_1,\ldots,n/\mu^n_I)=(\mu_1^{-1},\ldots,\mu_I^{-1}),\\\notag
\hat m&:=\limn \hat m^n=(\hat\lambda_i-\rho_i\hat\mu_i: i\in[I]),
\end{align}
and the $n\times n$ matrix
\begin{align}\notag
\sigma=(\sigma_{ij})_{1\le i,j\le I}:=\text{Diag}\left((\la_1)^{1/2} ,\ldots,(\la_I)^{1/2}\right).
\end{align}

We now define admissible controls for both players in the game, which due to their roles will be referred to as the {\it minimizer} and the {\it maximizer}.
\begin{itemize}
\item An {\it admissible control for the minimizer} is a tuple $\calE:=(\Omega,\calF,\{\calF_t\},\PP,\hat B,\hat Y)$, where $(\Omega,\calF,\{\calF_t\},\PP)$, 
is a filtered probability space
supporting the collection of independent one-dimensional $\calF_t$-adapted standard Brownian motions (SBMs) $\{\hat B_{j,i}: (j,i)\in\{A,S\}\times[I]\}$, and an $\R^I$-valued $\calF_t$-adapted process $\hat Y$ with RCLL sample paths such that $\theta\cdot \hat Y$ is nonnegative and nondecreasing.

\item
An {\it admissible control for the maximizer} is a measure $\Q$ defined on $(\Omega, \calF)$ such that for any $(j,i)\in\{A,S\}\times [I]$, 
\begin{align}\label{304}
\frac{d\Q_{j,i}}{d\PP_{j,i}}(t)=\exp\Big\{\int_0^t\hat \psi_{j,i}(s) d\hat B_{j,i}(s)-\frac{1}{2}\int_0^t \hat \psi^2_{j,i}(s)ds\Big\},\quad t\in\R_+,
\end{align}
for an $\calF_t$-progressively measurable process $\hat \psi_j:=(\hat\psi_{j,1},\ldots,\hat\psi_{j,I})$
satisfying
\begin{align}\label{305}
&\E^{\PP}\Big[\int_0^\iy \varrho(s)\hat \psi_{j,i}^2(s)ds\Big]<\iy\quad\text{and}\quad\E^{\PP}\Big[e^{\frac{1}{2}\int_0^t\hat \psi^2_{j,i}(s)ds}\Big]<\iy\quad t\in\supp(\varrho),\;i\in[I],
\end{align}
where $\Q_{j,i}:=\Q\circ(\hat B_{j,i})^{-1}$ and $\PP_{j,i}:=\PP\circ(\hat B_{j,i})^{-1}$. 
Moreover, the processes in the collection $\{\hat B_{j,i}: (j,i)\in\{A,S\}\times[I]\}$ are independent under $\Q$. 
\end{itemize}
Denote by $\calA$ the set of all admissible controls for the minimizer, where we often abuse notation and denote $\hat Y\in\calA$, keeping in mind that the control includes a filtered probability space. The set of all admissible controls for the maximizer is denoted by $\calQ$. Here as well we abuse terminology from time to time and refer to $\hat\psi$ as the maximizer's control.

Set $\hat B_j:=(\hat B_{j,1},\ldots, \hat B_{j,I})$, $j=A,S$, and let
\begin{align}\label{302}
\hat X(t)=\hat x_0+\hat mt+\sigma( \hat B_A(t)-\hat B_S(t)) +\hat Y(t),\quad t\in\R_+,
\end{align}
be the state process of the game
. An alternative form to the dynamics is as follows
\begin{align}
\notag
\hat X(t)=\hat x_0+\hat mt+\int_0^t\sigma[\hat\psi_A(s)-\hat\psi_S(s)]ds+\sigma(\hat B^{\Q}_A(t)-\hat B^{\Q}_S( t))+\hat Y(t),\quad t\in\R_+,
\end{align}
where $\hat B^{\Q}_j(\cdot):=\hat B_j(\cdot)-\int_0^\cdot\hat\psi_j(s)ds$, is an $I$-dimensional $\calF_t$-SBM under $\Q$.

%
%
%
Recall the definition of the cost in the QCP given in \eqref{212}. The cost associated with the strategy profile $(\hat Y,\Q)$ is given by
\begin{align}\notag
 J(\hat Y,\Q):= \;&
\E^{\Q}\Big[\int_0^\iy \varrho(t)C(\hat  X(t))dt \Big]-\sum_{i=1}^I\calL_{A,i}(\Q_{A,i}\|\PP_{A,i})-\sum_{i=1}^I\calL_{S,i}(\Q_{S,i}\|\PP_{S,i}).
\end{align}
The value function is thus
\begin{align}\notag
 V=\inf_{\hat Y\in\calA}\;\sup_{\Q\in\calQ}\; J(\hat Y,\Q)
\end{align}

\subsection{The solution of the game}\label{sec32}

In this section we provide a minimizing control for the minimizer in the game. For this, we present a key lemma regarding the minimizing curve.
\begin{lemma}\label{lem31}(Lemma 3.1 in \cite{ata-sah-cmu})
There exists a continuous function $f:\R_+\to\R_+^I$ such that for every $w\in\R_+$,
\begin{align}\label{308}
\theta\cdot f(w)=w \quad \text{and}\quad C(f(w))=\inf\{C(q):q\in\R_+^I,\theta\cdot q=w\}.
\end{align}
This function satisfies
\begin{align}\notag
\mu_1C_1'(f_1(w))=\ldots=\mu_IC_I'(f_I(w)).
\end{align}
Moreover, the mappings $w\mapsto C_i(f_i(w))$, $i\in[I]$, are increasing. 
\end{lemma}
One last ingredient for the definition of the candidate policy is the one-dimensional {\it Skorokhod map} $\Gamma:\calD(\R_+,\R)\to\calD(\R_+,\R)$ given by
\begin{align}\notag
\Gamma[l](t)=l(t)-\inf\{l(s)\wedge 0: s\in[0,t]\}, \quad t\in\R_+.
\end{align}
Pay attention that $\Gamma[l](t)\ge 0$ for every $t\in\R_+$. Moreover, it is well-known that for any $l_1,l_2\in\calD(\R_+,\R)$ and $t\in\R_+$,
\begin{align}\label{311}
|\Gamma[l_1]-\Gamma[l_2]|_t\le 2|l_1-l_2|_t.
\end{align}
Another nice feature of this function that serves us in the sequel is the {\it pathwise minimality property} of $\Gamma$. It says that for any $l,y\in \calD(\R_+,\R)$ such that $y$ is nonnegative and nondecreasing, and $l(t)+y(t)\ge0$ for all $t\in\R_+$, one has
\begin{align}\label{312}
l(t)+y(t)\ge\Gamma[l](t),\quad t\ge0.
\end{align}

Consider a filtered probability space as described in the previous subsection. 
Denote
\begin{align}
\notag
\hat L(t)=\hat x_0+\hat m t+\sigma  (\hat B_A(t)-\hat B_S(t))
\end{align}
 and set
\begin{align}\notag
\hat X_f(t)=f(\Gamma[\theta\cdot \hat L](t)))\quad \text{and} \quad\hat Y_f(t)=\hat X_f(t)-\hat L(t), \quad t\in\R_+.
\end{align}
We refer to the control $(\Omega,\calF,\{\calF_t\},\PP, \hat B, \hat Y_f)$ as the {\it $f$-reflecting control}. Occasionally, we abuse terminology and refer to $Y_f$ as the $f$-reflecting control.

\begin{proposition}\label{prop31}
The control $\calE:=(\Omega,\calF,\{\calF_t\},\PP, \hat B, \hat Y_f)$ is admissible and optimal for the minimizer. That is,
\begin{align}
\notag
V=\sup_{\Q\in\calQ}J(\hat Y_f,\Q).
\end{align}

\end{proposition}
The proof follows from the {\it pathwise minimality property} of the Skorokhod map combined with Lemma \ref{lem31} in the following sense. Let $\calE:=(\Omega',\calF',\{\calF_t'\}, \PP', \hat B', \hat Y')$ be an arbitrary admissible strategy for the minimizer and let $\hat X'$ be the associated dynamics via \eqref{302}. Since $\hat B=(\hat B_A,\hat B_S)$ and $\hat B'=((\hat B_A)',(\hat B_S)')$ are SBMs, we may couple between the two probability spaces such that both BMs are identified.

Now, by the definition of $f$, and \eqref{312}, for any $t\in\R_+$,
\begin{align}\notag
\theta\cdot \hat X_f(t)=\Gamma[\theta\cdot \hat L](t)\le \theta\cdot (\hat L(t)+\hat Y'(t))=\theta\cdot \hat X'(t). 
\end{align}
So, by \eqref{308}, 
 for any path of the BM $\omega$ and every $t\in\R_+$,
 \begin{align}\notag
C(\hat  X_f(t))(\omega)\le C(\hat X'(t))(\omega).
\end{align}
Hence, for any measure $\hat Q\in\calQ$,
\begin{align}\notag
 J( Y_f,\hat Q)\le  J( Y',\hat Q).
\end{align}

\begin{remark}
Note that the pathwise minimality property and the minimizing curve introduced in Lemma \ref{lem31} allow us to provide the optimal control for the minimizer in the limiting game without using differential equations and the dynamic programming principle. This is in contrast with \cite{coh2019a} where the differential equation is required to set the level of rejection and the best response of the maximizer player. The best response of the maximizer in the game plays an essential role in establishing the lower bound part of the convergence result in \cite{coh2019a}. Since our buffers are infinite, we are not required to the have a rejection level. Furthermore, as the last proof hints, the pathwise minimality property allows us to establish the lower bound in Section \ref{sec42} without a characterization of the optimal control for the maximizer in the limiting game. Hence, we do not perform a differential equations analysis.
\end{remark}

\section{Asymptotic optimality of the generalized $c\mu$-rule}\label{sec4}
\beginsec
\subsection{The generalized $c\mu$-rule and the main result}\label{sec41}
We now describe the generalized $c\mu$ rule. First we describe the preemptive version. This dynamic priority policy gives preemptive priority at time $t$ to
the class $i$ for which
$\mu_iC_i^{\prime}(\hat{X}^n_i(t))\geq \mu_jC_j^{\prime}(\hat{X}^n_j(t))$
for all $j$, where ties are broken in some arbitrary but predefined manner.
To define it precisely we need some additional notation.
Given a set of real numbers $R=\{r_i,i\in [I]\}$,
denote $\argmax R=\{i:r_i\ge\max_j r_j\}$, and
let $\argmax^*R$ be the smallest member of
$\argmax R$.
The control, that we denote by $U^{*,n}$, is defined by setting
\begin{align}\notag
U^{*,n}_i(t)=1_{\{\hat{X}^n(t)\in\calX_i\}}\,, \qquad i\in[I],
\end{align}
where $\calX_i$, $i\in[I]$ partition $\R_+^I\setminus\{0\}$ according to
\begin{equation}\notag
\calX_i=\{x\in\R_+^I\setminus\{0\}:\textstyle\argmax^*\{\mu_jC'_j(x_j),j\in[I]\}=i\},\qquad i\in[I].
\end{equation}
(Thus, in case of a tie, priority is given to the lowest index.)
Note that if, for some $i$, $x\in\calX_i$, we have $x_i>0$ thanks to the assumption that,
for all $i$, $C'_i(x)=0$ iff $x=0$; as a result, the queue selected for service is nonempty. \\
The non-preemptive version of the generalized $c\mu$ rule is a policy, denoted by
$U^{\#,n}$, that upon completion of a job selects a customer from the class $i$ for which
$\hat X^n\in\calX_i$. Namely, if $\tau$ is any time of departure, then $U^{\#,n}(\tau)=1_{\{\hat X^n(\tau)\in\calX_i\}}$.
Note that the job departing at time $\tau$ is not counted in $\hat X^n(\tau)$,
due to right-continuity. Both the policies have the non-idling property:
when a customer is admitted into an empty system, it
is immediately served. Now we state our main theorem.
 
\begin{theorem}
The value function of the QCP converges to the value function of the limiting game problem, i.e., $\lim_{n\to\iy}V^n=V$. Moreover, the generalized $c\mu$ rule (both preemptive as well as non-preemptive) is optimal. That is, if $\hat{Y}^n$ is the control corresponding to the generalized $c\mu$ rule (preemptive or non-preemptive), then 
$$\limsup_{n\to \iy}\sup_{\Q^n\in\calQ^n}  J^n(\hat Y^n,\Q^n)\leq V\,.$$
\end{theorem} 
The proof of the theorem is given in Sections \ref{sec42} and \ref{sec43}, where lower and upper bounds are established, respectively.

\subsection{Lower bound} \label{sec42}
By Proposition \ref{prop31} we need to show that
\begin{align}
\label{400}
\liminf_{n\to\iy} V^n\ge \sup_{\Q\in\calQ} J(\hat Y_f,\Q).
\end{align}
Or alternatively, for any arbitrary sequence of controls $\{\hat Y^n\}_n$,
\begin{align}
\label{401}
\liminf_{n\to\iy} \sup_{\Q^n\in\calQ^n} J^n(\hat Y^n,\Q^n)\ge \sup_{\Q\in\calQ} J(\hat Y_f,\Q).
\end{align}
Clearly, we may assume without loss of generality that the $\{\hat Y^n\}_n$ is taken s.t.~for every $n\in\N$,
\begin{align}\label{402}
 \sup_{\Q\in\calQ} J(\hat Y_f,\Q)+1\ge \sup_{\Q^n\in\calQ^n} J^n(\hat Y^n,\Q^n),
 \end{align}
otherwise, the lower bound holds trivially for such $\hat Y^n$. To this end, we need to follow these steps:
\begin{enumerate}
\item {\it Comparison between the two suprema.} Notice that the suprema on both sides of \eqref{402} are taken over very different sets. On the right-hand side (r.h.s.) the probability measures are w.r.t.~discrete processes and on the left-hand side (l.h.s.) the measures take care of the drift of the continuous process. In order to compare between the suprema 
we show that up to a small term $\eps>0$, the supremum on the l.h.s.~of \eqref{402} can be replaced by a supremum over measures whose corresponding $\psi$'s are nice functions of the Brownian motions. That is $\psi(t)=F^{\eps}(\hat B)$, for a nice function $F^{\eps}$. This way, we may compare between the two suprema by setting up a change of measure in the prelimit using the same function $F^{\eps}$, substituting the scaled arrival and departure processes in $F^{\eps}$ instead of $\hat B$.

\item {\it Tightness and convergence.} Showing tightness of a sequence of some relevant processes enables us to talk about converging sub-sequences. Restricting ourselves to such a subsequence, using the same function $F^{\eps}$ for the change of measure for the prelimit and limiting process, we show convergence of the divergence components. Moreover, using Lemma \ref{lem31} and the pathwise minimality property of the Skorokhod mapping \eqref{312} we bound from below the running cost of the prelimit problem, by the prelimit running cost associated with the $f$-reflecting control.
\end{enumerate}

The following lemma establishes the claim in the first part given above. Specifically, we show that for the limiting game there is an $\eps$-optimal $\hat \psi^\eps$ for the maximizer that is a bounded, and a continuous function of a finite sample of the BM, in a non-anticipating way. Its proof follows by the same arguments given in the proof of \cite[Theorem 10.3.1]{Kushner1992}, hence omitted.
\begin{lemma}\label{lem41}
For every $\eps>0$
	there is a system $\Xi^{\eps}\doteq (\Omega^{\eps},\calF^{\eps},\{\calF_t^{\eps}\},\PP^{\eps},\hat  B^{\eps})$
	and $\hat \psi^{\eps} \in \Pi(\Xi^\eps)$
	with the following properties.
	\begin{itemize}
		\item		$(\hat X^{\eps}_f,\hat Y^\eps_f)$ satisfies the following equation for every $t\in\R_+$.
\begin{align}\notag
\hat X^\eps_f(t)=f(\Gamma[\theta\cdot \hat L^\eps](t)))\quad \text{and} \quad\hat Y^\eps_f(t)=\hat X^\eps_f(t)-\hat L^\eps(t),
\end{align}
where $\hat L^\eps(t)=\hat x_0+\hat m t+\sigma( \hat B^{\eps}_A(t)-\hat B^{\eps}_S(t))$, and $\hat B^\eps=(\hat B^{\eps}_A,\hat B^{\eps}_S)$ is an $\calF_t$-adapted, $2I$-dimensional SBM under $\PP^\eps$.
		\item For some $\delta > 0$, $\hat \psi^{\eps}$ is piecewise constant on intervals of the form $[l\delta,(l+1)\delta)$, $l=0,1,2,\dots$. For every $s \in \R_+$, $\hat \psi^{\eps}(s)$ takes values in a finite subset of $\R^{2I}$, denoted by $Z^\eps$.
		\item For some $\theta>0$, for each $u \in Z^\eps$
\begin{align}\notag
\PP^{\eps}(\hat \psi^\eps(l\delta)=u\mid \hat B^\eps(s),s\le l\delta, \hat \psi^\eps(j\delta),j<l)&=\PP^{\eps}(\hat \psi^\eps(l\delta)=u\mid \hat B^\eps(p\theta),p\theta\le l\delta,\hat  \psi^\eps(j\delta),j<l)\nonumber\\
		&= F^{\eps}_u\left(\big(\hat B^\eps(p\theta)\big)_{p=0}^{\lfloor l\delta/\theta\rfloor} , \big(\hat \psi^\eps(j\delta)\big)_{j=0}^{l-1}\right),\notag
		\end{align}
		where 
		for suitable $ k_1, k_2 \in \N$, $F^{\eps}_u: \R^{ k_1} \times (Z^\eps)^{ k_2} \to [0,1]$ is a measurable function such that
		$F_u(\cdot, \mathbf{u})$ is continuous on $\R^{ k_1}$ for every $\mathbf{u} \in (Z^\eps)^{ k_2}$,
		\item Set the measure $\Q^\eps=\Pi_{i=1}^I(\Q^\eps_{A,i}\times\Q^\eps_{S,i})$ associated with $(\hat \psi^\eps_{j,i})_{j\in\{A,S\},i\in[I]}$ via \eqref{304}--\eqref{305}. Then,
		\begin{align}\label{405}
J(\hat Y^\eps_f,\Q^{\eps}) \ge \sup_{\Q\in\calQ}J(\hat Y^\eps_f,\Q) - \eps.
\end{align}
	\end{itemize}
\end{lemma}

Fix $\eps>0$ and $\Q^{\eps}$ such that \eqref{405} holds. The next proposition together with \eqref{401} establishes the lower bound \eqref{400}.
\begin{proposition} \label{prop41}
The following asymptotic bound holds
\begin{align}
\notag
\liminf_{n\to\iy}J(\hat Y^n,\Q^{n,\eps}) \ge J(\hat Y^\eps_f,\Q^\eps).
\end{align}
\end{proposition}
Before providing its proof, we need some notation and preleminary results.
	
For every $n\in\N$ set the process $(\hat \psi^{n,\eps}(t))$ to be random and fixed on the time interval $[l\delta,(l+1)\delta)$, that is a $Z^\eps$-valued, $\calF^n_t$-measurable according the the conditional distribution
\begin{align}\notag
\PP^n(\hat \psi^{n,\eps}(t)=u\mid\calF^n_t)=F^{\eps}_u\left(\hat M^n(p\theta)\big)_{p=0}^{\lfloor l\delta/\theta\rfloor}, \big(\hat \psi^{n,\eps}(j\delta)\big)_{j=0}^{l-1}\right),
\end{align}
where $\hat M^n:=(\hat A^n,\hat S^n(T^n))$. Let $\Q^{n,\eps}=\prod_{i=1}^I( \Q^{n,\eps}_{A,i}\times\Q^{n,\eps}_{S,i})$ be such that the measures $\Q^{n,\eps}_{A,i}$ and $\Q^{n,\eps}_{S,i}$ are respectively associated with the intensities $\psi^{n,\eps}_{A,i}$ and $\psi^{n,\eps}_{S,i}$, which are given by
\begin{align}\notag
\psi^{n,\eps}_{A,i}(t)&:=\la^n_i +\hat \psi^{n,\eps}_{A,i}(t)(\la_in)^{1/2},\qquad
\psi^{n,\eps}_{S,i}(t):=\mu^n_i +\hat \psi^{n,\eps}_{S,i}(t)(\mu_in)^{1/2},\qquad t\in\R_+.
\end{align}
Also, define the $\Q^{n,\eps}$-martingales
\begin{align}\label{408}
\begin{split}
\check A^n_i(t)&:=n^{-1/2}\Big(A^n_i(t)-\int_0^t \psi^{n,\eps}_{A,i}(s)ds\Big), \\
\check D^n_i(t)&:=n^{-1/2}\Big(S^n_i(T^n_i(t))-\int_0^t \psi^{n,\eps}_{S,i}(s)dT^n_i(s)\Big).
\end{split}
\end{align}
Pay attention that
\begin{align}
\notag 
\check A^n_i(t)=\hat A^n_i(t)-\la^{1/2}\int_0^t\hat \psi^{n,\eps}_{A,i}(s)ds,
\qquad\check D^n_i(t)=\hat S^n_i(T^n_i(t))-\mu^{1/2}\int_0^t\hat \psi^{n,\eps}_{S,i}(s)dT^n_i(s),
\end{align}
and therefore,
\begin{align}\label{410}
\hat X^n_i(t)=\hat X^n_i(0)+\hat m^n_i t+\check A^n_i(t)-\check D^n_i(t)+\hat Y^n_i(t)+\la_i^{1/2}\int_0^t\hat \psi^{n,\eps}_{A,i}(s)ds-\mu_i^{1/2}\int_0^t\hat \psi^{n,\eps}_{S,i}(s)dT^n_i(s).
\end{align}

The next lemma provides the first step in the proof of the second part.


\begin{lemma}\label{lem42}
For any given $\eps>0$, the sequence of processes $\{T^n\}_n$ converges in probability to $\brho$ under the measures $\{\Q^{n,\eps}\}_n$. 
\end{lemma}
{\bf Proof.} Throughout the proof we assume that $\varrho>0$. The case where $\varrho(t_0)=0$ for some $t_0>0$, and by monotonicity $\varrho(t)=0$ for every $t>t_0$ (finite horizon case), is treated similarly and therefore it is omitted. Set 
\begin{align}\notag
\theta^n:=\limn(n/\mu^n_1,\ldots,n/\mu^n_I). 
\end{align}
From \eqref{207} it is sufficient to show that for each $T>0$, the sequence $\{\Q^{n,\eps}\circ( \|\hat Y^n\|_T)^{-1}\}_n$ is tight.  This follows once we show that  the following two limits hold.
\begin{align}\label{411}
\lim_{K\to\iy}\limsup_{n\to\iy}\Q^{n,\eps}\left(\inf_{i\in[I]}\inf_{0\le t\le T}\hat Y^n_i(t)\le-K\right)&=0,\\\label{412}
\lim_{K\to\iy}\limsup_{n\to\iy}\Q^{n,\eps}\left(\theta^n\cdot\hat Y^n(T)\ge K\right)&=0.
\end{align}
To obeserve it, we now show that for sufficiently large $n$, the event $\|\hat Y^n\|_T\ge 2K/\theta_{\text{min}}$ implies that 
either (i) there exists $i\in[I]$ such that $\inf_{0\le t\le T}\theta^n_i\hat Y^n_i(t)\le -K/(4I)$ or (ii) $\theta^n\cdot\hat Y^n(T)\ge K/2$, where $\theta_{\text{min}}:=\min_{i\in[I]}\theta_i$.

Indeed, assume that $\|\hat Y^n\|_T\ge 2K/\theta_{\text{min}}$ 
and let $n$ be sufficiently large such that for every $i\in[I]$, $\theta^n_i\ge \theta_i-\theta_{\text{min}}/2$, where $\theta_{\text{min}}:=\min_{i\in[I]}\theta_i$. Now,
\begin{align}\notag
\sum_{i\in[I]}\|\theta_i^n\hat Y^n_i\|_T+\frac{\theta_{\text{min}}}{2}\sum_{i\in[I]}\|\hat Y^n_i\|_T\ge\sum_{i\in[I]}\|\theta_i\hat Y^n_i\|_T\ge\theta_{min}\sum_{i\in[I]}\|\hat Y^n_i\|_T.
\end{align}
Hence, 
\begin{align}\notag
\sum_{i\in[I]}\|\theta_i^n\hat Y^n_i\|_T\ge\frac{\theta_{\text{min}}}{2}\sum_{i\in[I]}\|\hat Y^n_i\|_T\ge\frac{\theta_{\text{min}}}{2}\|\hat Y^n\|_T.
\end{align}
Therefore, the event $\|\hat Y^n\|_T\ge 2K/\theta_{\text{min}}$ implies that $\sum_{i\in[I]}\|\theta^n_i\hat Y^n_i\|_T\ge K$. Now, either
\begin{enumerate}
\item[(1)] there exists $i\in[I]$ such that $\inf_{0\le t\le T}\theta^n_i\hat Y^n_i(t)\le -K/(4I)$; or
\item[(2)] for every $i\in[I]$, $\inf_{0\le t\le T}\theta^n_i\hat Y^n_i(t)> -K/(4I)$. This condition implies that $\theta^n\cdot\hat Y^n(T)\ge K/2$. Indeed, arguing by contradiction, assume that $\theta^n\cdot\hat Y^n(T)< K/2$. Now, $\theta^n\cdot\hat Y^n(T)=\sum_{i\in[I]}\theta^n_i\hat Y^n_i(T)$ can be decomposed into two partial sums, one that runs over the positive terms $\theta^n_i\hat Y^n_i(T)$ and the other over the negative ones. Denote them respectively by $\Sigma_+$ and $\Sigma_-$. Since $\theta^n\cdot\hat Y^n(T)=\Sigma_++\Sigma_-<K/2$, it follows that $\Sigma_+<-\Sigma_-+K/2< K/4+K/2$, where the last inequality follows by the case considered in this part. Therefore, $\sum_{i\in[I]}\|\theta^n_i\hat Y^n_i\|_T=\Sigma_+-\Sigma_-<K$, a contradiction to the conclusion mentioned above.
\end{enumerate}
Notice that for a sufficiently large $n$, for every $i\in[I]$, $\theta^n_i\le \theta_{\text{max}}+1$, where $\theta_{\text{max}}:=\max_{i\in[I]}\theta_i$. 
Thus, once \eqref{411}--\eqref{412} are established, we get,
\begin{align}\notag
&\limsup_{K\to\iy}\limsup_{n\to\iy}\Q^{n,\eps}\Big(\|\hat Y^n_i\|_T\ge2K/\theta_{\text{min}}\Big)\\\notag
&\quad\le
\limsup_{K\to\iy}\limsup_{n\to\iy}\Q^{n,\eps}\Big(\inf_{i\in[I]}\inf_{0\le t\le T}\hat Y^n_i(t)\le-K/(4I(\theta_{\text{max}}+1))\Big)
\\\notag
&\qquad+
\limsup_{K\to\iy}\limsup_{n\to\iy}\Q^{n,\eps}\Big(\theta^n\cdot\hat Y^n(T)\ge \frac{K}{2}\Big)\\\notag
&\quad=0.
\end{align}

\noi{\bf Establishing \eqref{411}:} 
Recall \eqref{410} and that $\hat X^n_i\ge 0$ and $\{\hat \psi^{n,\eps}_{j,i}\}_{j,i,n}$ are uniformly bounded. Hence, there exists a constant $a_1>0$ such that for every $n\in\N$ and $t\in[0,T]$,
\begin{align}\notag
\hat Y^n_i(t)\ge -a_1 -\check A^n_i(t)+\check D^n_i(t).
\end{align}
The event $\{\inf_{0\le t\le T}\hat Y^n_i(t)\le -K\}$ implies that 
either $\{\|\check A^n_i\|_T\ge (K-a_1)/2\}$ or $\{\|\check D^n_i\|_T\ge (K-a_1)/2\}$. By the Burkholder--Davis--Gundy (BDG) inequality there exists a constant $a_2>0$ such that 
\begin{align}\notag
\sup_n\E^{\Q^{n,\eps}}\left[\|\check A^n_i\|^2_T\right]\le
a_2\sup_n\Big\{n^{-1}\E^{\Q^{n,\eps}}\left[\| A^n_i(T)\|\right]\Big\}=:M_{A,i}<\iy,
\end{align}
where the last inequality follows since $A^n_i(T)$ is a Poisson random variable with mean $\E^{\Q^{n,\eps}}[\int_0^T\psi^{n,\eps}_{A,i}(t)dt]$ and since the processes $\{\hat \psi^{n,\eps}_{j,i}\}_{j,i,n}$ are uniformly bounded. A similar bound holds for $\check D^n_i$ with an associated constant $M_{S,i}<\infty$. Therefore,
\begin{align}
\notag
\Q^{n,\eps}\Big(\inf_{0\le t\le T}\inf_{i\in[I]}\hat Y^n_i(t)\le -K\Big)
&\le\sum_{i=1}^I\Big\{\Q^{n,\eps}\big(\|\check A^n_i\|_T\ge (K-a_1)/2\big)+\Q^{n,\eps}\big(\|\check D^n_i\|_T\ge (K-a_1)/2\big)\Big\}\\\notag
&\le \frac{2}{K-a_1}\sum_{i=1}^I(M_{A,i}+M_{S,i}).
\end{align}
The r.h.s.~converges to $0$ as $K\to\iy$.


\noi{\bf Establishing \eqref{412}:} 
Throughout this part, $a$ is a positive constant, independent of $n$ and $t$, which may change from one line to the next. Set $\tilde Y^n=(\tilde Y^n_i:i\in[I])$, $\bar Z^n=(\bar Z^n_i:i\in[I])$, 
$\bar A^n=(\bar A^n_i:i\in[I])$, and $\bar D^n=(\bar D^n_i:i\in[I])$ by
\begin{align}\notag
\tilde Y^n_i(t)=\frac{\theta^n_i}{\theta_i}\hat Y^n_i(t),\qquad\bar Z^n_i(t)= \Big|\frac{\theta^n_i}{\theta_i}-1\Big||\hat Y^n_i(t)|,\qquad\bar A^n_i(t)=|\check A^n_i(t)|,\qquad \bar D^n_i(t)=|\check D^n_i(t)|.
\end{align}
Also, set $e=(1,\ldots,1)\in\R^I$.
From \eqref{410} and the uniform boundedness of $\{\hat \psi^{n,\eps}_{j,i}\}_{j,i,n}$, we get that,
\begin{align}\notag
0\le \theta^n\cdot\hat Y^n(t)=\theta\cdot\tilde Y^n(t)
&\le \theta\cdot \left(\hat X^n(t)+ate+ \bar A^n(t)+ \bar D^n(t)+\bar Z^n(t)\right)\\\notag
&\le \theta\cdot \left(\hat X^n(t)+ate+ \bar A^n(t)+ \bar D^n(t)\right),
\end{align}
where the last inequality follows by modifying $a$, recalling that $|1-\theta^n_i/\theta_i|$ is of order $n^{-1/2}$ and that $|\hat Y^n_i(t)|\le a_3tn^{1/2}$ for some $a_3>0$ independent of $n,i$ and $t$.
Denote $\rho^{-1}=(\rho_i^{-1}:i\in[I])$. By the monotonicity of $C\circ f$, \eqref{308}, and the convexity of $C$, we get that
\begin{align}\label{418}
\begin{split}
&\E^{\Q^{n,\eps}}\Big[\int_0^{2T}\varrho(t)C(f(\tfrac{1}{4}\theta^n\cdot\hat Y^n(t)))dt\Big]
\\
&\quad\le
\E^{\Q^{n,\eps}}\Big[\int_0^{2T}\varrho(t)C(f(\theta\cdot(\tfrac{1}{4}[\hat X^n(t)+ate+ \bar A^n(t)+ \bar D^n(t)])))dt\Big]\\
&\quad\le
\E^{\Q^{n,\eps}}\Big[\int_0^{2T}\varrho(t)C(\tfrac{1}{4}[\hat X^n(t)+ate+ \bar A^n(t)+ \bar D^n(t)])dt\Big]\\
&\quad\le
\frac{1}{4}\left(\E^{\Q^{n,\eps}}\Big[\int_0^{2T}\varrho(t)C(\hat X^n(t))dt\Big]+
\E^{\Q^{n,\eps}}\Big[\int_0^{2T}\varrho(t)C(aet)dt\Big]\right.\\
&\qquad\qquad\left.+
\E^{\Q^{n,\eps}}\Big[\int_0^{2T}\varrho(t)C(\bar A^n(t)))dt\Big]+\E^{\Q^{n,\eps}}\Big[\int_0^{2T}\varrho(t)C(\bar D^n(t))dt\Big]\right).
\end{split}
\end{align}
Once we show that the r.h.s.~is uniformly bounded over $n$, we get by the monotonicity of $C_i\circ f_i$, $\theta^n\cdot \hat Y^n$, and $\varrho$, in addition to our assumption that $\varrho>0$, that
\begin{align}\notag
\sup_n
\E^{\Q^{n,\eps}}\Big[C_i(f_i(\tfrac{1}{4}\theta^n\cdot\hat Y^n(T)))\Big]
&\le a
\sup_n\E^{\Q^{n,\eps}}\Big[\int_{T}^{2T}\varrho(t)C_i(f_i(\tfrac{1}{4}\theta^n\cdot\hat Y^n(t)))dt\Big]
\\\notag
&\le a
\sup_n\E^{\Q^{n,\eps}}\Big[\int_0^{2T}\varrho(t)C_i(f_i(\tfrac{1}{4}\theta^n\cdot\hat Y^n(t)))dt\Big]<\iy.
\end{align}
By an application of Markov inequality and the monotonicity of $C_i\circ f_i$,
we obtain that \eqref{412} holds.

The rest of the proof is dedicated to uniformly bound the four terms on the r.h.s.~of \eqref{418}. The second term is deterministic and independent of $n$. Its bound follows by the polynomial growth of $C_i$ asserted in Assumption \ref{assumption1} and the last part of the assumption. To tackle the third term, we use again the polynomial growth of $C$ as follows
\begin{align}\notag
\sup_n\E^{\Q^{n,\eps}}\Big[\int_0^{2T}\varrho(t)C_i(|\check A^n_i(t)|)dt\Big]\le
a\sup_n\E^{\Q^{n,\eps}}\left[\|\check A^n_i\|^p_{2T}\right]\le
a\sup_nn^{-p/2}\E^{\Q^{n,\eps}}\left[\| A^n_i(2T)\|^{p/2}\right].
\end{align}
The last supremum is finite since $A^n_i(2T)$ is a Poisson random variable with mean $\int_0^{2T}\psi^{n,\eps}_{A,i}(t)dt$, the sequence $\{\hat \psi^{n,\eps}\}_n$ is uniformly bounded, and the $p/2$-moment of the Poisson random variables is a polynomial of order $p/2$.
The bound of the forth term is similar and therefore omitted.

%
%
%
%
%
%
%
%
%
In order to estimate the first expectation, recall \eqref{402}. Hence, for every $i\in[I]$ and $n\in\N$,,
\begin{align}\notag
\E^{\Q^{n,\eps}}\Big[\int_0^{2T}\varrho(t)C_i(\hat X^n_i(t))dt\Big]
\le
V+1+
 \sum_{i=1}^I\calL_{A,i}(\Q^{n,\eps}_{A,i}\|\PP^n_{A,i})
 +\sum_{i=1}^I\calL_{S,i}(\Q^{n,\eps}_{S,i}\|\PP^n_{S,i}).
 \end{align}
We now uniformly bound the divergence terms. We bound only the $\calL_{A,i}$-terms and the same arguments holds also for the $\calL_{S,i}$-terms. The growth condition of $g_{A,i}$ and simple algebraic manipulation of the Radon--Nykodim derivative from \eqref{216} yield that
\begin{align}\notag
&\E^{\Q^{n,\eps}}\Big[\int_0^{2T}\varrho(t)g_{A,i}\Big(\log\Big(\frac{d\Q^{n,\eps}_{A,i}}{d\PP^n_{A,i}}(t)\Big)\Big)dt\Big]\\\notag
&\quad\le
a+a\E^{\Q^{n,\eps}}\Big[\sup_{0\le t\le 2T}\Big|\log\Big(\frac{d\Q^{n,\eps}_{A,i}}{d\PP^n_{A,i}}(t)\Big)\Big|^{\bar{p}}\Big]\\\notag
&\quad\le
a+a\E^{\Q^{n,\eps}}\Big[\sup_{0\le t\le 2T}\Big|
\int_0^t\Big(\psi^{n,\eps}_{A,i}(s)\log\Big(\frac{\psi^{n,\eps}_{A,i}(s)}{\la^n_i}\Big) -\psi^{n,\eps}_{A,i}(s)+\la^n_i\Big)ds
\Big|^{\bar{p}}\Big]\\\notag
&\;\;\;\qquad
+a\E^{\Q^{n,\eps}}\Big[\sup_{0\le t\le 2T}\Big|\int_0^t\log\Big(\frac{\psi^{n,\eps}_{A,i}(s)}{\la^n_i}\Big) \Big(A^n_i(s)-\int_0^s\psi^{n,\eps}_{A,i}(u)du\Big)\Big|^{\bar{p}}\Big]
.
\end{align}
By the definition of $\psi^{n,\eps}_{A,i}$ and the inequality $(1+y)\log(1+y)-y\le y^2/2$ for $y$ in a neighbourhood of $0$ applied to $y^n=\hat \psi^{n,\eps}_{A,i}(t)(\la_in)^{1/2}/\la^n_i$ and the uniform bound over $\{\hat \psi^{n,\eps}_{A,i}\}_n$, one obtains that the first expectation on the r.h.s.~is bounded above. In order to estimate the second expectation pay attention that
\begin{align}\notag
\int_0^t\log\Big(\frac{\psi^{n,\eps}_{A,i}(s)}{\la^n_i}\Big)  d \Big(A^n_i(s)-\int_0^s\psi^{n,\eps}_{A,i}(u)du\Big)\quad\text{ is a martingale.}
\end{align}
Denote its quadratic variation, estimated at time $2T$ by $[\int_0^\cdot\dots]_{2T}$. Now, applying the BDG inequality and using the bound $|\log(1+y)|\le2|y|$ on a neighborhood of $0$ applied to $y^n$, we get that
\begin{align}\notag
&\E^{\Q^{n,\eps}}\Big[\sup_{0\le t\le 2T}\Big|\int_0^t\log\Big(\frac{\psi^{n,\eps}_{A,i}(s)}{\la^n_i}\Big)  d \Big(A^n_i(s)-\int_0^s\psi^{n,\eps}_{A,i}(u)du\Big)\Big|^{\bar{p}}\Big]
\\\notag
&\quad\le
a\E^{\Q^{n,\eps}}\left[\Big[\int_0^\cdot\log\Big(\frac{\psi^{n,\eps}_{A,i}(s)}{\la^n_i}\Big) d \Big(A^n_i(s)-\int_0^s\psi^{n,\eps}_{A,i}(u)du\Big)\Big]^{{\bar{p}}/2}_{2T}\right]
\\\notag
&\quad\le
an^{-{\bar{p}}/2}\E^{\Q^{n,\eps}}[(A^n_i(2T))^{{\bar{p}}/2}]<\iy
.
\end{align}
The last bound follows since $A^n_i(2T)$ is a Poisson random variable with mean $\int_0^{2T}\psi^{n,\eps}_{A,i}(t)dt$, the sequence $\{\hat \psi^{n,\eps}_{A,i}\}_n$ is uniformly bounded, and the ${\bar{p}}/2$-moment of the Poisson random variables is a polynomial of order ${\bar{p}}/2$.

\qed

The lower bound will be established via weak convergence and tightness arguments. For this, we set up the rest of the processes required for this purpose. We start with breaking the logarithm of the Radon--Nykodim derivatives \eqref{216}, \eqref{216z}, and \eqref{304} into two parts each. For every $n\in\N$ and $i\in[I]$, set the processes $H^n_{A,i}, G^n_{A,i}, H^n_{S,i},$ and $G^n_{S,i}$ by
\begin{align}\notag
H^n_{A,i}(t)&=\int_0^tn^{1/2}\log\Big(\frac{\psi^{n,\eps}_{A,i}(s)}{\la^n_i}\Big)d\check A^n(s),\\\notag
G^n_{A,i}(t)&=\int_0^t\Big(\psi^{n,\eps}_{A,i}(s)\log\Big(\frac{\psi^{n,\eps}_{A,i}(s)}{\la^n_i}\Big) -\psi^{n,\eps}_{A,i}(s)+\la^n_i\Big)ds,\\\notag
H^n_{S,i}(t)&=\int_0^tn^{1/2}\log\Big(\frac{\psi^{n,\eps}_{S,i}(s)}{\mu^n_i}\Big)d\check D^n(s),\\\notag
G^n_{S,i}(t)&=\int_0^t\Big(\psi^{n,\eps}_{S,i}(s)\log\Big(\frac{\psi^{n,\eps}_{S,i}(s)}{\mu^n_i}\Big) -\psi^{n,\eps}_{S,i}(s)+\mu^n_i\Big)dT^n_i(s).
\end{align}
Moreover, set the processes $H_{A,i}, G_{A,i}, H_{S,i},$ and $G_{S,i}$ by
\begin{align}\notag
H_{A,i}(t)&=\int_0^t\hat \psi^{\eps}_{A,i}(s)d\hat B^{\eps,\Q^{\eps}}_{A,i}(s),\qquad G_{A,i}(t)=\frac{1}{2}\int_0^t(\hat \psi^{\eps}_{A,i}(s))^2ds,\\\notag
H_{S,i}(t)&=\int_0^t\hat \psi^{\eps}_{S,i}(s)d\hat B^{\eps,\Q^{\eps}}_{S,i}(s),\qquad G_{S,i}(t)=\frac{1}{2}\int_0^t(\hat \psi^{\eps}_{S,i}(s))^2ds,
\end{align}
where
$$(\hat B^{\eps,\Q^{\eps}}_{A,i},\hat B^{\eps,\Q^{\eps}}_{S,i})(\cdot):=(B^{\eps}_{A,i},B^{\eps}_{S,i})(\cdot)-\int_0^\cdot(\hat \psi^\eps_{A,i}(s),\hat \psi^\eps_{S,i}(s))ds.$$
Denote $H^n=(H^n_{j,i}:j=A,S; i\in[I]), G^n=(G^n_{j,i}:j=A,S; i\in[I])$ and similarly $H=(H_{j,i}:j=A,S;t i\in[I]), G=(G_{j,i}:j=A,S; i\in[I])$. 
Furthermore, recall that we aim at bounding the limit inferior of $V^n$ by the robust cost associated with generalized $c\mu$ control, which in turn is defined via reflection. To reach this cost, we need to pass through a process with a reflection structure in the prelimit. Hence, we set up the following.  Let $\hat L^n=(\hat L^n_i:i\in[I]), \hat X^n_f,$ and $\hat Y^n_f$ be defined as follows
\begin{align}\notag
\hat L^n_i(t)=\hat X^n_i(0)+\hat m^n_i t+\hat A^n_i(t)-\hat S^n_i(T^n_i(t))+n^{-1/2}(\mu^n_i - n\mu_i)(\rho_i t-T^n_i(t)),
\end{align}
%
and
\begin{align}\notag
\hat X^n_f(t)=f(\Gamma[\theta\cdot\hat L^n(\cdot)](t)),\qquad
\hat Y^n_f(t)=\hat X^n_f(t)-\hat L^n(t).
\end{align}
\begin{lemma} \label{lem43}
The following sequence of measures
\begin{align}\label{429}
\left\{\Q^{n,\eps}\circ\left(\hat A^n,\hat S^n, T^n,\hat S^n(T^n),\check A^n,\check D^n,\{\hat \psi^{n,\eps}_{j,i}\}_{j,i}, H^n, G^n,\hat L^n, \hat X^n_f, \hat Y^n_f\right)^{-1}\right\}_n
\end{align}
is $\calC$-tight. Moreover, every limit point of this sequence has the same distribution as
\begin{align}\label{430}
\Q^\eps\circ\left(\sigma\hat B^{\eps}_A,(\rho)^{-1}\sigma\hat B^\eps_S,\brho,\sigma\hat B^\eps_S, \sigma\hat B^{\eps,\Q^{\eps}}_A,\sigma\hat B^{\eps,\Q^{\eps}}_S,\{\hat \psi^\eps_{j,i}\}_{j,i},H,G,\hat L^\eps,\hat X^\eps_f,\hat Y^\eps_f\right)^{-1},
\end{align}
where, the process $(\hat B^{\eps,\Q^{\eps}}_{A,i},\hat B^{\eps,\Q^{\eps}}_{S,i})$ is a $2I$-dimensional SBM under the measure $\Q^\eps$ and the filtration $\calF$ that is generated by the processes in \eqref{430}.
\end{lemma}
{\bf Proof.} The tightness argument is standard and therefore omitted. As for the limit, using the martingale central limit theorem as well as Lemma \ref{lem41} we obtain the convergence of all the terms besides
that of $(\hat \psi^{n,\eps}, H^n, G^n)$. To show the convergence of the latter, notice that the continuity of $F^\eps_u$ implies $\Q^{n,\eps}\circ(\hat \psi^{n,\eps})^{-1}\To \Q^\eps\circ(\hat \psi^\eps)^{-1}$. By the definition of $\psi^{n,\eps}$, the uniformly boundedness of $\{\hat \psi^{n,\eps}\}_n$, and the martingale central limit theorem, we finally obtain that $\Q^\eps\circ(H^n,G^n)^{-1}\To\Q^\eps\circ(H,G)^{-1}$.

\qed

{\bf Proof of Proposition \ref{prop41}.}
From the previous lemma we may reduce to a converging subsequence of \eqref{429}, which we relabel by $\{n\}$. In order to establish the desired lower bound it is sufficient to prove the following asymptotic estimates:  
\begin{align}\label{431}
&\limn \Big\{\sum_{i=1}^I\calL_{A,i}(\Q^{n,\eps}_{A,i}\|\PP^n_{A,i})
 +\sum_{i=1}^I\calL_{S,i}(\Q^{n,\eps}_{S,i}\|\PP^n_{S,i})\Big\}\\\notag
 &\qquad\qquad\qquad=
 \sum_{i=1}^I\calL_{A,i}(\Q_{A,i}\|\PP^\eps_{A,i})-\sum_{i=1}^I\calL_{j,i}(\Q_{S,i}\|\PP^\eps_{S,i})
\end{align}
and
\begin{align}\label{432}
&\liminf_{n\to\iy}\int_0^\iy \varrho(t)C(\hat X^n(t))dt\ge\int_0^\iy \varrho(t)C(\hat X_f(t))dt.
\end{align}

The limit \eqref{431} follows from Lemma \ref{lem43} and the representations \begin{align}\notag
 \calL_{j,i}(\Q^{n,\eps}_{j,i}\|\PP^n_{j,i})
 &=\E^{\Q^{n,\eps}}\Big[\int_0^\iy \varrho(t)g_{j,i}(H^n_{j,i}(t)+G^n_{j,i}(t))dt\Big],
\\\notag
 \calL_{j,i}(\Q^{\eps}_{j,i}\|\PP^\eps_{j,i})
 &=\E^{\Q^{\eps}}\Big[\int_0^\iy \varrho(t)g_{j,i}(H_{j,i}(t)+G_{j,i}(t))dt\Big].
\end{align}

We now turn to proving the lower bound \eqref{432}. The idea in this part is to use the properties of $f$ given in Lemma \ref{lem31} and the pathwise minimality property of the Skorokhod map, see \eqref{312}, applied to $\theta\cdot\hat L^n$. For this, notice that
\begin{align}\notag
\theta\cdot\hat X^n(t)=\theta\cdot\hat L^n(t)+n^{1/2}\Big(1-\sum_{i=1}^IT^n_i(t)\Big),
\end{align}
with $\theta\cdot \hat X^n\ge 0$ and $\Big(1-\sum_{i=1}^IT^n_i(t)\Big)$ is nonnegative and nondecreasing. Hence, \eqref{312} implies that
\begin{align}\notag
\theta\cdot \hat X^n(t)\ge \Gamma[\theta\cdot\hat L^n(\cdot)](t),\quad t\in[0,T].
\end{align}
From Lemma \ref{lem31}, the monotonicity of $C\circ f$, and the last bound, we get that
\begin{align}\notag
&\int_0^\iy \varrho(t)C(\hat X^n(t))dt\\\notag
&\quad
\ge \int_0^\iy \varrho(t)C(f(\theta\cdot\hat X^n(t)))dt\\\notag
&\quad
\ge\int_0^\iy \varrho(t)C(f(\Gamma[\theta\cdot\hat L^n(\cdot)](t))dt\\\notag
&\quad
= \int_0^\iy \varrho(t)C(\hat X^n_f(t))dt.
\end{align}
Finally, Lemma \ref{lem43} implies that
\begin{align}\notag
&\liminf_{n\to\iy}\int_0^\iy \varrho(t)C(\hat X^n(t))dt
\ge \int_0^\iy \varrho(t)C(\hat X^\eps_f(t))dt.
\end{align}

%
%
%
%
  \qed

\subsection{Upper bound}\label{sec43}

In this part we show that the generalized $c\mu$ rule asymptotically attains the value of the limiting problem $V$.
We denote by $\hat Y^n$ the control corresponding to the generalized $c\mu$ rule (preemptive or non-preemptive) and show that
\begin{align}\label{4381}
\limsup_{n\to \iy}\sup_{\Q^n\in\calQ^n}  J^n(\hat Y^n,\Q^n)\leq V\;.
\end{align}
For this we set up in Section \ref{sec431} an arbitrary sequence of measures $\{\Q^n\}_n$ and show that for any ``reasonable" sequence from the point of view of the maximizer, the $p$-means of the drifts and the logarithms of the Radon--Nykodym derivatives of the $n$-th systems are uniformly bounded. Then, in Section \ref{sec432} we use this uniform bound to show that the measures $\{\Q^n\}_n$ can be uniformly approximated by measures $\{\Q^{n,k}\}_n$ for some sufficiently large $k>0$, such that the associated rates $\psi^{n,k}_{j,i}$ satisfy,
$$\psi^{n,k}_{A,i}(t)=\la^n_i+(\la_in)^{1/2} \hat\psi^{n,k}_{A,i}(t)+o(n^{-1/2}),$$
where $|\psi^{n,k}_{A,i}|\le k$, and similarly for $\psi^{n,k}_{S,i}$. The motivation behind this step is that while $\int_0^\cdot\hat\psi^n_{j,i}(t)dt$ has a convergence subsequence, the limit is not necessarily absolutely continuous, hence, the limit might not have an integral form $\int_0^\cdot\hat\psi_{j,i}(s)ds$. Hence, we cannot compare it with the limiting game.  Moreover, we use it to obtain the convergence of the Radon--Nykodim derivatives. 
Finally, in Section \ref{sec433} we characterize the limiting process by providing the state-space collapse and in Section \ref{sec434} we establish the upper bound.

\subsubsection{$p$-mean bound for the intensities}
\label{sec431}
Consider an arbitrary sequence of measures chosen by the maximizer in the QCP, $\{\Q^n\}_n$, satisfying \eqref{216}--\eqref{216z} for some $\{\psi^n_{j,i}\}_{n,j,i}$. The first step in establishing the truncation reduction is showing that the maximizer can be restricted to measures $\Q^n\in\calQ^n$, which are close in average to the reference measure $\PP^n$. The reason it holds is because for high values of $n^{-1/2}(|\psi^n_{A,i}-\lambda^n_i|+|\psi^n_{S,i}-\mu^n_i|)$, the divergence terms (in absolute values) significantly dominate the running costs, which are affected through $\hat X^n$. Without loss of generality, we may assume that for any $n\in\N$,
\begin{align}\label{439}
J^n(\hat Y^n,\Q^n)\ge V-1.
\end{align}
Otherwise, the upper bound holds trivially for such $\Q^n$'s. 

Recall the definitions of $\check A^n$ and $\check D^n$ given in \eqref{408} and set
\begin{align}
\notag
\hat\psi^n_{A,i}(t):=(\la_i n)^{-1/2}\left(\psi^n_{1,i}(t)-\la^n_i\right),\qquad\text{and}\qquad\hat\psi^n_{S,i}(t):=(\mu_i n)^{-1/2}\left(\psi^n_{2,i}(t)-\mu^n_i\right).
\end{align}
Notice that
\begin{align}\label{440}
\begin{split}
&\calL_{A,i}(\Q^n_{A,i}\|\PP^n_{A,i})\\
&\quad=\E^{\Q^n}\Big[\int_0^\iy \varrho (t)g_{A,i}\Big(\int_0^tf^n_{A,i}(\hat\psi^n_{A,i}(s))ds+\int_0^tn^{1/2}\log\Big(\frac{\psi^n_{A,i}(s)}{\la^n_i}\Big)d\check A^n_i(s)\Big)dt\Big],\\
&\calL_{S,i}(\Q^n_{S,i}\|\PP^n_{S,i})\\
&\quad=
\E^{\Q^n}\Big[\int_0^\iy \varrho(t)g_{A,i}\Big(\int_0^tf^n_{S,i}(\hat\psi^n_{S,i}(s))dT^n_i(s)+\int_0^tn^{1/2}\log\Big(\frac{\psi^n_{S,i}(s)}{\mu^n_i}\Big)d\check D^n_i(s)\Big)dt\Big],
\end{split}
\end{align}
where
$f^n_{A,i}:\left(-\la^n_i(\la_in)^{-1/2},\iy\right)\to\R_+$ and $f^n_{S,i}:\left(-\mu^n_i(\mu_in)^{-1/2},\iy\right)\to\R_+$ are given by
\begin{align}
\notag
\begin{split}
f^n_{A,i}(x):=\la^n_i\Big[\left(1+\tfrac{(\la_in)^{1/2}}{\la^n_i}x\right)\log\left(1+\tfrac{(\la_in)^{1/2}}{\la^n_i}x)\right)
-\tfrac{(\la_in)^{1/2}}{\la^n_i}x\Big],\\
f^n_{S,i}(x):=\mu^n_i\Big[\left(1+\tfrac{(\mu_in)^{1/2}}{\mu^n_i}x\right)\log\left(1+\tfrac{(\mu_in)^{1/2}}{\mu^n_i}x)\right)
-\tfrac{(\mu_in)^{1/2}}{\mu^n_i}x\Big].
\end{split}
\end{align}
In the sequel, we need the following properties of $f^n_{j,i}$, $j\in\{A,S\},i\in[I]$:
\begin{itemize}
\item the function $x\mapsto f^n_{j,i}(x)/x$ is increasing on $(0,\iy)$,
\item $\lim_{x\to\iy}\sup_n f^n_{j,i}(x)/x=\infty$.
\end{itemize}

Pay attention that we do not assume that the processes $\sup_{j,i, n}|\hat\psi^n_{j,i}(t)|$ are uniformly bounded. Rather, we use the next proposition to claim that one may restrict these processes to be uniformly bounded without too much loss.
\begin{proposition}\label{prop42}
There exists $M>0$ 
such that for every $n\in\N$ and every
$i\in[I]$,
\begin{align}\label{444}
&\E^{\Q^n}\Big[\int_0^\iy\varrho(t)\Big\{\Big|\int_0^t\hat\psi^n_{A,i}(s)ds\Big|^{\bar p}+\Big|\int_0^t\hat\psi^n_{S,i}(s)dT^n_{S,i}(s)\Big|^{\bar p}\Big\}dt\Big]
\le M,
\end{align}
\begin{align}\label{445}
&E^{\Q^n}\Big[\int_0^\iy \varrho (t)\Big\{\Big(\int_0^tf^n_{A,i}(\hat\psi^n_{A,i}(s))ds\Big)^{\bar p}+\Big(\int_0^tf^n_{S,i}(\hat\psi^n_{S,i}(s))dT^n_i(s)\Big)^{\bar p}\Big\}dt\Big]
\le M,
\end{align}
and \begin{align}\label{445b}
&E^{\Q^n}\Big[\int_0^\iy \varrho (t)(\hat X^n_i(t))^{\bar p}dt\Big]
\le M.
\end{align}

\end{proposition}
{\bf Proof.} Throughout the proof, the parameter $a$ stands for a positive constant that is independent of $n$ and $t$ and which can change from one line to the next. Pay attention that $0<n^{-1}\psi^n_{j,i}(t)\le a(1+n^{-1/2}\hat\psi^n_{j,i}(t))$, $t\in\R_+$.
Applying the BDG inequality to $\check A^n$ and $\check D^n$, we have
\begin{align}
\label{446}
\begin{split}&\E^{\Q^n}[\|\check A^n_i\|_t^p]\le an^{-p/2}\E^{\Q^n}\Big[\Big|\int_0^t\psi^n_{A,i}(s)ds\Big|^{p/2}\Big]\le a(t^{p/2}+n^{-p/4}\E^{\Q^n}\Big[\Big|\int_0^t\hat \psi^n_{A,i}(s)ds\Big|^{p/2}\Big]\Big),\\
&\E^{\Q^n}[\|\check D^n_i\|_t^p]\le an^{-p/2}\E^{\Q^n}\Big[\Big|\int_0^t\psi^n_{S,i}(s)dT^n_i(s)\Big|^{p/2}\Big]\le a\Big(t^{p/2}+n^{-p/4}\E^{\Q^n}\Big[\Big|\int_0^t\hat \psi^n_{S,i}(s)dT^n_i(s)\Big|^{p/2}\Big]\Big).
\end{split}
\end{align}
Recall that both versions of the generalized $c\mu$-policy, preemptive and non-preemptive, are work conserving, that is $\sum_{i\in[I]}U^{*,n}_i(t)=1$ and $\sum_{i\in[I]}U^{\sharp,n}_i(t)=1$ whenever $\hat X^{n,k}(t)$ is nonzero. By the definitions of $T^{n,k}$ and $\hat Y^{n,k}$ it follows that the nondecreasing process $\theta^n\cdot Y^{n,k}$, does not increase when $\theta^n\cdot X^{n,k}>0$. 
From \eqref{410} we get that for any $t\in\R_+$, 
\begin{align}\label{447}
&\theta^n\cdot\hat X^n(t)\\\notag
&\quad=\Gamma\Big[\theta^n\cdot\Big(\hat X^n(0)+\hat m^n\cdot+\check A^n(\cdot)+\check D^n(\cdot)+\int_0^\cdot\sigma\hat\psi^n_{A,i}(s)ds-\int_0^\cdot\sigma_S\hat\psi^n_{S,i}(s)dT^n_i(s)\Big)\Big](t),
\end{align}
where $\sigma_S:=\text{Diag}(\mu^{1/2}_1,\ldots,\mu^{1/2}_I)$. By \eqref{311}, the above, the uniform bound $0\le \hat X^n_i(t)\le a\theta^n\cdot\hat X^n(t)$, and since $\{\theta^n\}_n$ is uniformly bounded, it follows that for any $t\in\R_+$,
%
%
\begin{align}\label{tmp1}
(\hat X^n_i(t))^p
&\le a\sum_{i=1}^I\Big(1+t^p+(\|\check A^n_i\|_t)^p+(\|\check D^n_i\|_t)^p
+\Big|\int_0^t\hat\psi^n_{A,i}(s)ds\Big|^p+\Big|\int_0^t\hat\psi^n_{S,i}(s)dT^n_{S,i}(s)\Big|^p
\Big).
\end{align}

By the polynomial growth of the running cost and the bound $\int_0^\iy\varrho(t)t^pdt<\iy$, both given in Assumption \ref{assumption1}, 
\begin{align}\notag
&\sum_{i=1}^I\E^{\Q^n}\Big[\int_0^\iy\varrho(t)C_i(\hat X^n_i(t))dt\Big]
\\\notag
&\quad\le
a\sum_{i=1}^I\Big(1+\E^{\Q^n}\Big[\int_0^\iy\varrho(t)\Big\{\Big|\int_0^t\hat\psi^n_{A,i}(s)ds\Big|^p+\Big|\int_0^t\hat\psi^n_{S,i}(s)dT^n_{S,i}(s)\Big|^p\Big\}dt\Big]\Big).
\end{align}
Combining it with the definition of the cost function $J^n$ and \eqref{439} one obtains that
\begin{align}\label{450}
&\calL_{A,i}(\Q^n_{A,i}\|\PP^n_{A,i})+\calL_{S,i}(\Q^n_{S,i}\|\PP^n_{S,i})\\\notag
&\qquad
\le
a\Big(1+\E^{\Q^n}\Big[\int_0^\iy\varrho(t)\Big\{\Big|\int_0^t\hat\psi^n_{A,i}(s)ds\Big|^p+\Big|\int_0^t\hat\psi^n_{S,i}(s)dT^n_{S,i}(s)\Big|^p\Big\}dt\Big]\Big).
\end{align}
Pay attention that for any $t\in\R_+$,
\begin{align}\notag
\E^{\Q^n}\Big[\int_0^tn^{1/2}\log\Big(\frac{\psi^n_{A,i}(s)}{\la^n_i}\Big)d\check A^n_i(s)\;\big|\; \calF^{\psi^n}_t\Big]=0, 
\end{align}
where $\calF^{\psi^n}_t:=\sigma\{\psi^n_{j,i}(s):s\le t\;, j=A,S\;,i\in[I]\}$.
Also, recall that the functions $g_{j,i}$, $j\in\{A,S\}, i\in[I]$ are convex. Hence, by Jensen's inequality
\begin{align*}
&\E^{\Q^n}\Big[g_{A,i}\Big(\int_0^tf^n_{A,i}(\hat\psi^n_{A,i}(s))ds+\int_0^tn^{1/2}\log\Big(\frac{\psi^n_{A,i}(s)}{\la^n_i}\Big)d\check A^n_i(s)\Big)\Big]\\\notag
&\quad\ge\E^{\Q^n}\Big[g_{A,i}\Big(\E^{\Q^n}\Big\{\int_0^tf^n_{A,i}(\hat\psi^n_{A,i}(s))ds+\int_0^tn^{1/2}\log\Big(\frac{\psi^n_{A,i}(s)}{\la^n_i}\Big)d\check A^n_i(s)\;\big|\; \calF^{\psi^n}_t\Big\}\Big)\Big]\\
&\quad=\E^{\Q^n}\Big[g_{A,i}\Big(\E^{\Q^n}\Big\{\int_0^tf^n_{A,i}(\hat\psi^n_{A,i}(s))ds\;\big|\;\calF^{\psi^n}_t\Big\}\Big)\Big]\\
&\quad=\E^{\Q^n}\Big[g_{A,i}\Big(\int_0^tf^n_{A,i}(\hat\psi^n_{A,i}(s))ds\Big)\Big],
\end{align*}
and similarly by conditioning on $\sigma\{\psi^n_{j,i}(s), T^n_i(s):s\le t\;, j=A,S\;,i\in[I]\}$,
\begin{align*}
&\E^{\Q^n}\Big[g_{S,i}\Big(\int_0^tf^n_{S,i}(\hat\psi^n_{S,i}(s))dT^n_i(s)+\int_0^tn^{1/2}\log\Big(\frac{\psi^n_{S,i}(s)}{\mu^n_i}\Big)d\check D^n_i(s)\Big)\Big]\\\notag
&\quad\ge\E^{\Q^n}\Big[g_{S,i}\Big(\int_0^tf^n_{S,i}(\hat\psi^n_{S,i}(s))dT^n_i(s)\Big)\Big].
\end{align*}
Plugging in the expressions of the divergences given in \eqref{440} and using the two bounds above together with the bound $g_{j,i}(x)\ge c_1+c_2x^{\bar p}$ that Assumption \ref{assumption1} asserts, \eqref{450} yields that
\begin{align}\label{453}
&E^{\Q^n_{A,i}}\Big[\int_0^\iy \varrho (t)\Big(\int_0^tf^n_{A,i}(\hat\psi^n_{A,i}(s))ds\Big)^{\bar p}dt\Big]+
\E^{\Q^n_{S,i}}\Big[\int_0^\iy \varrho(t)\Big(\int_0^tf^n_{S,i}(\hat\psi^n_{S,i}(s))dT^n_i(s)\Big)^{\bar p}dt\Big]\\\notag
&\qquad
\le
a\Big(1+\E^{\Q^n}\Big[\int_0^\iy\varrho(t)\Big\{\Big|\int_0^t\hat\psi^n_{A,i}(s)ds\Big|^p+\Big|\int_0^t\hat\psi^n_{S,i}(s)dT^n_{S,i}(s)\Big|^p\Big\}dt\Big]\Big).
\end{align}

Denote for every $t\in\R_+$,
\begin{align}\notag
y^n_{A,i}(t):=\frac{(\la_in)^{1/2}}{\la^n_i}\hat\psi^n_{A,i}(t),\qquad y^n_{S,i}(t):=\frac{(\mu_in)^{1/2}}{\mu^n_i}\hat\psi^n_{S,i}(t).
\end{align}
Now, one can simply verify that for all $y> -1$,
\begin{align}\notag
\frac{1}{4}y^2\one_{\{y<4\}}+y\one_{\{y\ge4\}}\le(1+y)\log(1+y)-y.
\end{align}
Using this inequality and the one from \eqref{453} in addition to the definitions of $\la^n_i$ and $\mu^n_i$ given in \eqref{205}, we get that there exists $a_1>0$, such that for any $n\in\N$ and $t\in\R_+$,
\begin{align}\notag
&\E^{\Q^n}\Big[\int_0^\iy \varrho (t)\Big(\int_0^t\Big\{(\hat\psi^n_{A,i}(s))^2\one_{\{y^n_{A,i}(s)<4\}}+n^{1/2}\hat\psi^n_{A,i}(s)\one_{\{y^n_{A,i}(s)\ge4\}}\Big\}ds\Big)^{\bar p}dt\Big]
\\\notag
&\quad
+\E^{\Q^n}\Big[\int_0^\iy\varrho (t)\Big(\int_0^t\Big\{(\hat\psi^n_{S,i}(s))^2\one_{\{y^n_{S,i}(s)<4\}}
+n^{1/2}\hat\psi^n_{S,i}(s)\one_{\{y^n_{S,i}(s)\ge4\}}
\Big\}dT^n_i(s)\Big)^{\bar p}dt\Big]\\\notag
&\qquad
\le a_1\Big(1+\E^{\Q^n}\Big[\int_0^\iy\varrho(t)\Big\{\Big|\int_0^t\hat\psi^n_{A,i}(s)ds\Big|^p+\Big|\int_0^t\hat\psi^n_{S,i}(s)dT^n_{S,i}(s)\Big|^p\Big\}dt\Big]\Big).
\end{align}
Since the mapping $\psi\mapsto \psi^2$ is super-linear, there is a constant $a_2<0$ such that for any $\psi\in\R$, $\psi^2\ge a_2+2a_1|\psi|$. Applying this inequality for $\hat\psi^n_{j,i}(t)$ on the left-hand side of the above, we get that for every $n\ge 4a_1^2$,
\begin{align}\notag
\begin{split}
&a_2+2a_1\Big(\E^{\Q^n}\Big[\int_0^\iy\varrho(t)\Big\{\Big|\int_0^t\hat\psi^n_{A,i}(s)ds\Big|^{\bar p}+\Big|\int_0^t\hat\psi^n_{S,i}(s)dT^n_{S,i}(s)\Big|^{\bar p}\Big\}dt\Big]\Big)\\
&\quad\le
 a_1\Big(1+\E^{\Q^n}\Big[\int_0^\iy\varrho(t)\Big\{\Big|\int_0^t\hat\psi^n_{A,i}(s)ds\Big|^p+\Big|\int_0^t\hat\psi^n_{S,i}(s)dT^n_{S,i}(s)\Big|^p\Big\}dt\Big]\Big)
.
\end{split}
\end{align}
Recall also that $\bar p\ge p$, then,
\begin{align}\notag
\E^{\Q^n}\Big[\int_0^\iy\varrho(t)\Big\{\Big|\int_0^t\hat\psi^n_{A,i}(s)ds\Big|^{\bar p}+\Big|\int_0^t\hat\psi^n_{S,i}(s)dT^n_{S,i}(s)\Big|^{\bar p}\Big\}dt\Big]\le 1-a_2/a_1,
\end{align}
and \eqref{444} is established. The bound in \eqref{445} follows by another application of \eqref{453} and the bound above. Finally, the bound in \eqref{445b} follows by combining the bounds from \eqref{444}, \eqref{446}, $\int_0^\iy\varrho(t)t^{\bar p}dt<\iy$,  and \eqref{tmp1}.

\qed

\subsubsection{Reduction to uniformly truncated intensities}\label{sec432}

Having at hand the uniform bound (over the expectations) from the previous proposition, we now claim that up to a small loss from the maximizer's point of view , the terms $\{\hat\psi^n_{j,i}\}_{j,i,n}$ can be uniformly bounded. For this, we set up for every $k>0$ the processes
\begin{align}\notag
\psi^{n,k}_{A,i}(t)&:=\psi^n_{A,i}(t)-(\la_in)^{1/2}\hat\psi^n_{A,i}(t)\one_{\{|\hat\psi^n_{A,i}(t)|>k\}},\\\notag
\psi^{n,k}_{S,i}(t)&:=\psi^n_{S,i}(t)-(\mu_in)^{1/2}\hat\psi^n_{S,i}(t)\one_{\{|\hat\psi^n_{S,i}(t)|>k\}}.
\end{align}
Also, denote by $T^{n,k}=(T^{n,k}_i:i\in[I])$  the DM's generalized $c\mu$ rule given in Section \ref{sec41} associated with the environment associated with the intensities $\{\psi^{n,k}_{j,i}\}_{j,i}$. The arrival and service processes associated with these truncated intensities are $A^{n,k}$ and $S^{n,k}$, and they are coupled with $A^n$ and $S^n$ as follows. 
Set the following independent Poisson processes (with rate 1): $\{P_{j,i,m}:j\in\{A,S\},i\in[I], m=1,\ldots, 4\}$. 
For every $i\in[I]$ set up the following processes
\begin{align}\notag
&M^{n,-}_{A,i}(\cdot)=P_{A,i,1}\Big(\int_0^\cdot \Big(\la^n +n^{1/2}\hat\psi^n_{A,i}(s)\one_{\{\hat\psi^n_{A,i}(s)<0\}}\Big)ds\Big),\\\notag
&K^{n,-}_{A,i}(\cdot)=P_{A,i,2}\Big(\int_0^\cdot n^{1/2}\Big(-\hat\psi^n_{A,i}(s))\one_{\{\hat\psi^n_{A,i}(s)<-k\}}\Big)ds\Big),\\\notag
&K^{n,+}_{A,i}(\cdot)=P_{A,i,3}\Big(\int_0^\cdot n^{1/2}\hat\psi^n_{A,i}(s)\one_{\{0<\hat\psi^n_{A,i}(s)\le k\}}ds\Big),\\\notag
&M^{n,+}_{A,i}(\cdot)=P_{A,i,4}\Big(\int_0^\cdot n^{1/2}\hat\psi^n_{A,i}(s)\one_{\{\hat\psi^n_{A,i}(s)>k\}}ds\Big),
\end{align}
and similarly,
\begin{align}\notag
&M^{n,-}_{S,i}(\cdot)=P_{S,i,1}\Big(\int_0^\cdot (\mu^n +n^{1/2}\hat\psi^n_{S,i}(s)\one_{\{\hat\psi^n_{S,i}(s)<0\}})dT^n_i(s)\Big),\\\notag
&M^{n,k,-}_{S,i}(\cdot)=P_{S,i,1}\Big(\int_0^\cdot (\mu^n +n^{1/2}\hat\psi^n_{S,i}(s)\one_{\{\hat\psi^n_{S,i}(s)<0\}})dT^{n,k}_i(s)\Big),\\\notag
&K^{n,k,-}_{S,i}(\cdot)=P_{S,i,2}\Big(\int_0^\cdot n^{1/2}(-\hat\psi^n_{S,i}(s))\one_{\{\hat\psi^n_{S,i}(s)<-k\}}dT^{n,k}_i(s)\Big),\\\notag
&K^{n,+}_{S,i}(\cdot)=P_{S,i,3}\Big(\int_0^\cdot n^{1/2}\hat\psi^n_{S,i}(s)\one_{\{0<\hat\psi^n_{S,i}(s)\le k\}}dT^n_i(s)\Big),\\\notag
&K^{n,k,+}_{S,i}(\cdot)=P_{S,i,3}\Big(\int_0^\cdot n^{1/2}\hat\psi^n_{S,i}(s)\one_{\{0<\hat\psi^n_{S,i}(s)\le k\}}dT^{n,k}_i(s)\Big),\\\notag
&M^{n,+}_{S,i}(\cdot)=P_{S,i,4}\Big(\int_0^\cdot n^{1/2}\hat\psi^n_{S,i}(s)\one_{\{\hat\psi^n_{S,i}(s)>k\}})dT^n_i(s)\Big).
\end{align}
Now set, $A^n=(A^n_i:i\in[I])$, $A^{n,k}=(A^{n,k}_i:i\in[I])$, $D^n=(D^n_i:i\in[I])$, $D^{n,k}=(D^{n,k}_i:i\in[I])$ as follows
\begin{align}\notag
A^n&:=M^{n,-}_A+K^{n,+}_A+M^{n,+}_A,\qquad A^{n,k}:=M^{n,-}_A+K^{n,-}_A+K^{n,+}_A\\\notag
D^n&:=M^{n,-}_S+K^{n,+}_S+M^{n,+}_S,\qquad D^{n,k}:=M^{n,k,-}_S+K^{n,k,-}_S+K^{n,k,+}_S.
\end{align}
Also, denote by $\check A^n, \check A^{n,k},\check D^n$, and $\check D^{n,k}$  the compensated versions of $A^n, A^{n,k},D^n$, and $D^{n,k}$, respectively.
Finally, set $\hat Y^{n,k}=(\hat Y^{n,k}_i:i\in[I])$ with $\hat Y^{n,k}(\cdot)=\mu^n_in^{-1/2}(\rho_i\cdot-T^{n,k}_i(\cdot))$ and the state process
\begin{align}\notag
\hat X^{n,k}_i(t)&=\hat X^n_i(0)+\hat m^n_i t+\check A^{n,k}_i(t)-\check D^{n,k}_i(t)+\hat Y^{n,k}_i(t)\\\notag
&\quad+\la_i^{1/2}\int_0^t\hat \psi^{n,k}_{A,i}(s)ds-\mu_i^{1/2}\int_0^t\hat \psi^{n,k}_{S,i}(s)dT^{n,k}_i(s).
\end{align}
As in \eqref{447}, we have as well for any $t\in\R_+$,
\begin{align}\label{463}
&\theta^n\cdot\hat X^{n,k}(t)\\\notag
&\quad=\Gamma\Big[\theta^n\cdot\Big(\hat X^n(0)+\hat m^n\cdot+\check A^{n,k}(\cdot)+\check D^{n,k}(\cdot)+\int_0^\cdot\sigma\hat\psi^{n,k}_{A,i}(s)ds-\int_0^\cdot\sigma_S\hat\psi^{n,k}_{S,i}(s)dT^{n,k}_i(s)\Big)\Big](t),
\end{align}
where recall that $\sigma_S=\text{Diag}(\mu^{1/2}_1,\ldots,\mu^{1/2}_I)$.

\begin{lemma} \label{lem44}
For any given $k>0$, the sequence of processes $\{(T^n,T^{n,k})\}_n$ converges in probability to  $(\brho,\brho)$ under the measures $\{\Q^{n}\}_n$. 
\end{lemma}
Once we establish the convergence in probability for each of the components, the joint convergence follows. As in the proof of Lemma \ref{lem42}, in order to establish the convergence of each of the components, it is sufficient to prove that for every $T>0$, $\{\Q^n\circ((\|\hat Y^n\|_T)^{-1})$ is tight. The proof here is similar, where now the uniform boundedness of $\int_0^T\hat\psi^{n,\eps}_{A,i}(t)dt$, asserted in Section \ref{sec42},  is replaced by the uniform boundedness of $\E^{\Q^n}\big[\big|\int_0^T\hat\psi^n_{A,i}(t)dt\big|^{\bar p}\big]$ and similarly for $j=S$. The proof is therefore, omitted.
\begin{proposition}\label{prop43}
The following asymptotic bound holds
\begin{align}\notag
\lim_{k\to\iy}\;
\liminf_{n\to\iy} \;\big\{J(\hat Y^{n,k},\Q^{n,k})-J(\hat Y^{n},\Q^{n})\big\}\ge 0.
\end{align}
\end{proposition}

\noi{\bf Proof.} Throughout the proof, the parameter $a$ stands for a positive constant, independent of $n,t$, and $k$, and which can change from one line to the next. The proof is done in two parts, separately taking care of the holding costs and the divergence terms.

\skp\noi{\it Part (i)}.
We start with showing that 
\begin{align}\label{465}
\liminf_{n\to\iy}\E^{\Q^n}\Big[\int_0^\iy\varrho(t)\big\{C(\hat X^{n,k}(t))-C(\hat X^n(t))\big\}dt\Big]\ge 0.
\end{align}
The convexity of $C_i$ implies that for every $i\in[I]$,
\begin{align}\notag
&\E^{\Q^n}\Big[\int_0^\iy\varrho(t)\big\{C_i(\hat X^{n,k}_i(t))-C_i(\hat X^n_i(t))\big\}dt\Big]\\\notag
&\quad \ge 
\E^{\Q^n}\Big[\int_0^\iy\varrho(t)C_i'(\hat X^n_i(t))(\hat X^{n,k}_i(t)-\hat X^n_i(t))dt\Big].
\end{align}

Assumption \ref{assumption1} implies that $C'(x)\le a(1+x^{p-1})$. 
Hence, it is sufficient to show that for every $i\in[I]$,
\begin{align}\notag
\lim_{n\to\iy}\E^{\Q^n}\Big[\int_0^\iy\varrho(t)\left(1+(\hat X^n_i(t))^{p-1}\right)\big|\hat X^{n,k}_i(t)-\hat X^n_i(t)\big|dt\Big]= 0.
\end{align}
By H{\"o}lder's inequality (using the powers $p/(p-1)$ and $p$) and \eqref{445b} it is sufficient to show that 
\begin{align}\notag
\lim_{n\to\iy}\E^{\Q^n}\Big[\int_0^\iy\varrho(t)\big|\hat X^{n,k}_i(t)-\hat X^n_i(t)\big|^pdt\Big]= 0.
\end{align}
Moreover, since, $\{\theta^n_i\}_{i,n}$ are bounded away from $0$, the latter follows once we show that 
\begin{align}\label{466}
\lim_{n\to\iy}\E^{\Q^n}\Big[\int_0^\iy\varrho(t)|\theta^n\cdot\hat X^n(t)-\theta^n\cdot\hat X^{n,k}(t)|^pdt\Big]=0.
\end{align}
Pay attention that by the BDG inequality,
\begin{align}\notag
\E^{\Q^n}
\Big[
\|\check M^{n,-}-\check M^{n,k,-}\|^p_t
&\le 
 an^{-p/2}\E^{\Q^n}\Big[\Big\{
\int_0^t\Big(\mu^n +n^{1/2}\hat\psi^n_{S,i}(s)\one_{\{\hat\psi^n_{S,i}(s)<0\}}\Big)d|T^n_i(s)-dT^{n,k}_i(s)|\Big\}^{p/2}
\Big]\\\notag
&\le
a\E^{\Q^n}
\Big[
\|T^n_i-T^{n,k}_i\|^{p/2}_t
\Big].
\end{align}
By \eqref{447} and \eqref{463},
\begin{align}
\notag
&\E^{\Q^n}\left[\big|\theta^n\cdot\hat X^n(t)-\theta^n\cdot\hat X^{n,k}(t)\big|^p\right]\\\notag
&\quad\le
a\sum_{i=1}^I\Big\{\E^{\Q^n}\Big[\Big(\int_0^t|\hat\psi^n_{A,i}(s)|\one_{\{|\hat\psi^n_{A,i}(s)|>k\}}ds\Big)^p\Big]
\\\notag
&\quad\qquad\qquad+\E^{\Q^n}\Big[\Big(\int_0^t|\hat\psi^n_{S,i}(s)|\one_{\{|\hat\psi^n_{S,i}(s)|>k\}}dT^n_i(s)\Big)^p\Big]
\\\notag
&\quad\qquad\qquad+\E^{\Q^n}\Big[\Big(\int_0^t\big(|\hat\psi^{n,k}_{S,i}(s)|\big)d\big(T^n_i(s)-T^{n,k}_i(s)\big)\Big)^p\Big]
\\\notag
&\quad\qquad\qquad
+\E^{\Q^n}
\Big[
\|T^n_i-T^{n,k}_i\|^{p/2}_t\Big\}
.
\end{align}
Therefore,
\begin{align}\label{468}
\begin{split}
&\E^{\Q^n}\Big[\int_0^\iy\varrho(t)|\theta^n\cdot\hat X^n(t)-\theta^n\cdot\hat X^{n,k}(t)|^pdt\Big]]\\
&\quad\le a\sum_{i=1}^I\Big\{\E^{\Q^n}\Big[\int_0^\iy\varrho(t)\Big(\int_0^t|\hat\psi^n_{A,i}(s)|\one_{\{|\hat\psi^n_{A,i}(s)|>k\}}ds\Big)^pdt\Big]
\\
&\qquad\qquad\quad+
\E^{\Q^n}\Big[\int_0^\iy\varrho(t)\Big(\int_0^t|\hat\psi^n_{S,i}(s)|\one_{\{|\hat\psi^n_{S,i}(s)|>k\}}dT^n_i(s)\Big)^pdt\Big]
\\
&\qquad\qquad\quad+
\E^{\Q^n}\Big[\int_0^\iy\varrho(t)\Big(\int_0^t|\hat\psi^{n,k}_{S,i}(s)|d\big(T^n_i(s)-T^{n,k}_i(s)\big)\Big)^pdt\Big]
\\
&\qquad\qquad\quad+
\E^{\Q^n}\Big[\int_0^\iy\varrho(t)\|T^n_i-T^{n,k}_i\|^{p/2}_tdt\Big]\Big\}.
\end{split}
\end{align}
We show that by taking $\lim_{k\to\iy}\liminf_{n\to\iy}$ the four terms of the sum on the r.h.s.~of the above converge to zero. The convergence of the first two terms follow by the same argument, which for convenience, we provide only for the first one. For every $t\in\R_+$,
\begin{align}\notag
\begin{split}
&\int_0^t|\hat\psi^n_{A,i}(s)|\one_{\{|\hat\psi^n_{A,i}(s)|>k\}}ds
\\
&\quad=
\int_0^t\hat\psi^n_{A,i}(s)\one_{\{\hat\psi^n_{A,i}(s)>k\}}ds
+\int_0^t-\hat\psi^n_{A,i}(s)\one_{\{-\hat\psi^n_{A,i}(s)>k\}}ds
\\
&\quad
\le\frac{k}{f^n_{A,i}(k)}\int_0^tf^n_{A,i}(\hat\psi^n_{A,i}(s))\one_{\{\hat\psi^n_{A,i}(s)>k\}}ds
+
a\frac{n}{k}\int_0^t\Big(\frac{(\la_in)^{1/2}\hat\psi^n_{A,i}(s)}{\la^n_i}\Big)^2\one_{\{\hat\psi^n_{A,i}(s)<-k\}}ds
\\
&\quad
\le\Big(\frac{k}{f^n_{A,i}(k)}+\frac{a}{k}\Big)\int_0^tf^n_{A,i}(\hat\psi^n_{A,i}(s))ds.
\end{split}
\end{align}
The first inequality follows since the function $x\mapsto f^n_{j,i}(x)/x$ is increasing, since for $x>k$, $x\le x^2/k$, and since $(\la_i n)^{1/2}/\la^n_i$ is or order $n^{-1/2}$. The second one follows since for $x<0$, $x^2\le (1+x)\log(1+x)-x$ and again since $(\la_i n)^{1/2}/\la^n_i$ is of order $n^{-1/2}$. Taking $\E^{\Q^n}[\int_0^\iy\varrho(t)(\dotsm)^{\bar p} dt]$ on both sides and using \eqref{445} and the limit $\lim_{k\to\iy} \sup_n f^n_{j,i}(k)/k=\iy$, one obtains the convergence of the first term of \eqref{468}.
%
%
%
For any given $k>0$, the third term on the r.h.s.~of \eqref{468} converges to zero as $n\to\iy$ since $|\hat \psi^{n,k}_{A,i}|\le k$ and by Lemma \ref{lem44}. Finally, the last term on the r.h.s.~of \eqref{468} converges to zero by Lemma \ref{lem44}. These limits imply that \eqref{466} holds, which in turn implies that \eqref{465} holds.

\skp\noi {\it Part (ii)}. 
We now turn to the divergence terms. 
We show that for any $i\in[I]$,
\begin{align}\notag
&\lim_{k\to\iy}\limsup_{n\to\iy}\Big\{
\calL_{A,i}(\Q^{n,k}_{A,i}\|\PP^n_{A,i})
-
\calL_{A,i}(\Q^n_{A,i}\|\PP^n_{A,i})
\Big\}\le 0,\\\notag
&\lim_{k\to\iy}\limsup_{n\to\iy}\Big\{
\calL_{S,i}(\Q^{n,k}_{S,i}\|\PP^n_{S,i})
-
\calL_{S,i}(\Q^n_{S,i}\|\PP^n_{S,i})
\Big\}\le 0.
\end{align}
The proofs of both asymptotic bounds are similar, however, the components $T^n$ and $T^{n,k}$ add another level of complication to the proof of the second limit. Hence, we only prove the latter.
Denote
\begin{align}\notag
E^n(t)&:=\int_0^tf^n_{S,i}(\hat\psi^n_{S,i}(s))dT^n_i(s),\qquad F^n(t):=\int_0^tn^{1/2}\log\Big(\frac{\psi^n_{S,i}(s)}{\mu^n_i}\Big)d\check D^n_i(s),\\\notag
E^{n,k}(t)&:=\int_0^tf^n_{S,i}(\hat\psi^{n,k}_{S,i}(s))dT^{n,k}_i(s),\qquad F^{n,k}(t):=\int_0^tn^{1/2}\log\Big(\frac{\psi^{n,k}_{S,i}(s)}{\mu^n_i}\Big)d\check D^{n,k}_i(s).
\end{align}

The conditional expectation of $F^n(t)$ given $\calF^{\psi^n}_t=\sigma\{\psi^n_{j,i}(s):s\le t\;, j=A,S\;,i\in[I]\}$ is zero. 
From \eqref{440} and the convexity of $g_{S,i}$ if follows that
\begin{align}\label{473}
&\calL_{S,i}(\Q^n_{S,i}\|\PP^n_{S,i})-\calL_{S,i}(\Q^{n,k}_{S,i}\|\PP^n_{S,i})\\\notag
&=
\E^{\Q^n}\Big[\E^{\Q^n}\Big\{\int_0^\iy \varrho(t)g_{S,i}\left(E^n(t)+F^n(t)\right)dt\;\Big|\;\calF^{\psi^n}\Big\}\Big]\\\notag
&\quad
-\E^{\Q^n}\Big[\int_0^\iy \varrho(t)\big\{g_{S,i}\left(E^{n,k}(t)+F^{n,k}(t)\right)dt\Big]\\\notag
&\ge
\E^{\Q^n}\Big[\int_0^\iy \varrho(t)g_{S,i}\left(\E^{\Q^n}\big\{E^n(t)+F^n(t)\mid\calF^{\psi^n}_t\big\}\right)dt\Big]
\\\notag
&\quad
-\E^{\Q^n}\Big[\int_0^\iy \varrho(t)\big\{g_{S,i}\left(E^{n,k}(t)+F^{n,k}(t)\right)dt\Big]\\\notag
&\ge\E^{\Q^n}\Big[\int_0^\iy \varrho(t)g_{S,i}'\left(E^{n,k}(t)+F^{n,k}(t)\right)\left(\E^{\Q^n}\big\{E^n(t)\mid\calF^{\psi^n}_t\big\}-E^{n,k}(t)-F^{n,k}(t)\right)dt\Big]
.
\end{align}
We now show that $\lim_{k\to\iy}\liminf_{n\to\iy}$ of the r.h.s.~is 0. 
To this end, we show that
\begin{align}\label{lim1}
&\lim_{k\to\iy}\liminf_{n\to\iy}\E^{\Q^n}\Big[\int_0^\iy \varrho(t)g_{S,i}'\left(E^{n,k}(t)+F^{n,k}(t)\right)F^{n,k}(t)dt\Big]=0,\\\label{lim2}
&\lim_{k\to\iy}\liminf_{n\to\iy}\E^{\Q^n}\Big[\int_0^\iy \varrho(t)g_{S,i}'\left(E^{n,k}(t)+F^{n,k}(t)\right)\left(\E^{\Q^n}\big\{E^n(t)\mid\calF^{\psi^n}_t\big\}-E^{n,k}(t)\right)dt\Big]\ge 0.
\end{align}
 
We start with the proof of \eqref{lim1}. The conditions on $g_{S,i}$ imply that its derivative has at most a polynomial growth of order $\bar p-1$. Applying H\"older's inequality as in the first part of the proof, we get that it is sufficient to show that 
\begin{align}\notag
&\limsup_{k\to\iy}\limsup_{n\to\iy} \E^{\Q^n}\Big[\int_0^\iy \varrho(t)|E^{n}(t)|^{\bar p}dt\Big]<\iy,
\end{align}
and\begin{align}\label{4731}
&\lim_{k\to\iy}\limsup_{n\to\iy}\E^{\Q^n}\Big[\int_0^\iy \varrho(t)|F^{n,k}(t)|^{\bar p}dt\Big]=0.
\end{align}
The first bound follows from Proposition \ref{prop42}. To obtain the second bound we first apply the BDG inequality and obtain that
\begin{align}\notag
\E^{\Q^n}\Big[\int_0^\iy \varrho(t)|F^{n,k}(t)|^{\bar p}dt\Big]\le a\E^{\Q^n}\Big[\Big(\int_0^tn\Big|\log\Big(\frac{\psi^n_{S,i}(s)}{\mu^n_i}\Big)\Big|^2\one_{\{|\hat\psi^n_{S,i}(s)|\le k\}}d \tfrac{D^{n,k}_i(s)}{n}\Big)^{\bar p/2}\Big].
\end{align}
Pay attention that for any given $k>0$, there is $n_k>0$ such that for all $n\ge n_k$, and any $s\in\R_+$, 
\begin{align}\notag
n^{1/2}\Big|\log\Big(\frac{\psi^n_{S,i}(s)}{\mu^n_i}\Big)\Big|\one_{\{|\hat\psi^n_{S,i}(s)|\le k\}}\le 2|\hat\psi^n_{S,i}(s)|.
\end{align}
Therefore, the r.h.s.~of the above is bounded above by
\begin{align}\notag
\E^{\Q^n}\Big[\int_0^\iy \varrho(t)|F^{n,k}(t)|^{\bar p}dt\Big]\le a\E^{\Q^n}\Big[\Big(\int_0^t\big|\hat\psi^n_{S,i}(s)\big|^2\one_{\{|\hat\psi^n_{S,i}(s)|\le k\}}d \tfrac{D^{n,k}_i(s)}{n}\Big)^{\bar p/2}\Big].
\end{align}
Finally, Jensen's inequality together with the order $n$ rate of $D^{n,k}$ and the same arguments given at the end of part (i) of this proof, relying on Proposition \ref{prop42}, imply that \eqref{4731} holds and \eqref{lim1} is established.

We now turn to the proof of \eqref{lim2}. 
Pay attention that 
\begin{align}\notag
\E^{\Q^n}\big[E^n(t)\mid\calF^{\psi^n}_t\big]-E^{n,k}(t) 
&=
\E^{\Q^n}\Big[\int_0^tf^n_{S,i}(\hat\psi^n_{S,i}(s))\one_{\{|\hat\psi^n_{S,i}(t)|> k\}}dT^n_i(s)\;\Big|\;\calF^{\psi^n}_t\Big]\\\notag
&\quad
+
\E^{\Q^n}\Big[\int_0^tf^n_{S,i}(\hat\psi^{n,k}_{S,i}(s))d(T^n_i(s)-\rho_is)\;\Big|\calF^{\psi^n}_t\Big]\\\notag
&\quad
+
\int_0^tf^n_{S,i}(\hat\psi^{n,k}_{S,i}(s))d(\rho_is-T^{n,k}_i(s))
.
\end{align}
Since $g_{j,i}'$ is non-decreasing we can use the last display together with the fact that $f^n_{j,i}$ is nonnegative to get
\begin{align}\notag
&\lim_{k\to\iy}\liminf_{n\to\iy}\E^{\Q^n}\Big[\int_0^\iy \varrho(t)g_{S,i}'\left(E^{n,k}(t)+F^{n,k}(t)\right)\left(\E^{\Q^n}\big[E^n(t)\mid\calF^{\psi^n}_t\big]-E^{n,k}(t) \right)dt\Big]\\\notag
&\quad\ge
\lim_{k\to\iy}\liminf_{n\to\iy}\E^{\Q^n}
\Big[\int_0^t\Big\{\varrho(t)g_{S,i}'\left(E^{n,k}(t)+F^{n,k}(t)\right)\E^{\Q^n}\Big[\int_0^tf^n_{S,i}(\hat\psi^{n,k}_{S,i}(s))d(T^n_i(s)-\rho_is)\;\Big|\calF^{\psi^n}_t\Big]\Big\}dt\Big]\\\notag
&\qquad
+\lim_{k\to\iy}\liminf_{n\to\iy}\E^{\Q^n}
\Big[\int_0^t\Big\{\varrho(t)g_{S,i}'\left(E^{n,k}(t)+F^{n,k}(t)\right)\int_0^tf^n_{S,i}(\hat\psi^{n,k}_{S,i}(s))d(\rho_is-T^{n,k}_i(s))\Big\}dt\Big].
\end{align}
We now show that the last two limits are zero. Since the proofs for both limits are similar, we focus only on the second one and show that
\begin{align}\notag
\lim_{k\to\iy}\limsup_{n\to\iy}\E^{\Q^n}
\Big[\int_0^t\Big\{\varrho(t)g_{S,i}'\left(E^{n,k}(t)+F^{n,k}(t)\right)B^{n,k}(t)\Big\}dt\Big]=0,
\end{align}
where 
$$B^{n,k}(t):=\Big|\int_0^tf^n_{S,i}(\hat\psi^{n,k}_{S,i}(s))d(\rho_is-T^{n,k}_i(s))\Big|.$$
Now, 
The conditions on $g_{S,i}$ imply that its derivative has at most a polynomial growth of order $\bar p-1$. Again applying H\"older's inequality, we get that it is sufficient to show that 
\begin{align}\notag
&\limsup_{k\to\iy}\limsup_{n\to\iy} \E^{\Q^n}\Big[\int_0^\iy \varrho(t)|E^{n,k}(t)|^{\bar p}dt\Big]<\iy,\\\notag
&\limsup_{k\to\iy}\limsup_{n\to\iy} \E^{\Q^n}\Big[\int_0^\iy \varrho(t)|F^{n,k}(t)|^{\bar p}dt\Big]<\iy,
\end{align}
and\begin{align}\label{4731}
&\lim_{k\to\iy}\limsup_{n\to\iy}\E^{\Q^n}\Big[\int_0^\iy \varrho(t)|B^{n,k}(t)|^{\bar p}dt\Big]=0.
\end{align}
The first two bounds were established before. Finally, the last limit follows by Lemma \ref{lem44} and since for any given $k$, $\sup_nf^n_{S,i}(\hat\psi^{n,k}_{S,i}(s))$ is uniformly bounded.


\qed

\subsubsection{State-space collapse}\label{sec433}
We now focus on the truncated processes. Fix as arbitrary $k>0$ and define 
\begin{align}\notag
\hat W^{n,k}(t):=\check A^{n,k}(t)-\check D^{n,k}(t)+\hat m^nt,\qquad 
\hat\Psi^{n,k}_j=\int_0^t\hat\psi^{n,k}_j(s)ds, \qquad, j\in\{A,S\},\quad t\in\R_+.
\end{align}

We now state the state-space collapse property. Its proof follows by the same arguments given in \cite{ata-sah-cmu, ASAFCOHEN2018, van} and relies on the fact that on any given compact time interval $\Psi^{n,k}_j$, $j=A,S$, are uniformly bounded (for the fixed $k$). Therefore, it is omitted.
\begin{proposition}\label{prop44}
The following limit holds
\begin{align}\notag
\limn\Q^{n,k}\circ\left(\hat X^{n,k}-f(\theta^n\cdot \hat X^{n,k})\right)^{-1}=0.
\end{align}
\end{proposition}
Next we provide a limiting result.
\begin{lemma}\label{lem45}
The following sequence of probability measures is $\calC$-tight 
\begin{align}\notag
&\left\{\Q^{n,k}\circ
\left(
\check A^{n,k},\check D^{n,k},\hat W^{n,k},\hat X^{n,k}, \hat Y^{n,k}, \hat\Psi^{n,k}_A,\hat\Psi^{n,k}_S, 
T^{n,k},\Big\{\frac{d\Q^{n,k}_{j,i}}{d\PP^n_{j,i}}\Big\}_{i,j}
\right)^{-1},\right.
\\\notag
&\quad\left.
\PP^{n,k}\circ
\left(
\check A^{n,k},\check D^{n,k},\hat W^{n,k},\hat X^{n,k}, \hat Y^{n,k}, \hat\Psi^{n,k}_A,\hat\Psi^{n,k}_S, 
T^{n,k},\Big\{\frac{d\Q^{n,k}_{j,i}}{d\PP^n_{j,i}}\Big\}_{i,j}
\right)^{-1}
\right\}_n
\end{align}
and any sub-sequential limit of it
\begin{align}\notag
&\Q^{\circ,k}\circ
\left(\check A^{\circ,k},\check D^{\circ,k},
\hat W^{\circ,k},\hat X^{\circ,k}, \hat Y^{\circ,k}, \hat\Psi^{\circ,k}_A,\hat\Psi^{\circ,k}_S, 
\brho,\{H_{j,i}\}_{j,i}
\right)^{-1},\\\notag
&\quad
\PP^{\circ,k}\circ
\left(\check A^{\circ,k},\check D^{\circ,k},
\hat W^{\circ,k},\hat X^{\circ,k}, \hat Y^{\circ,k}, \hat\Psi^{\circ,k}_A,\hat\Psi^{\circ,k}_S, 
\brho,\{H_{j,i}\}_{j,i}\Big\}_{j,i}
\right)^{-1}
\end{align}
satisfies 
\begin{enumerate}
\item $\hat X^{\circ,k}(0)=\hat X(0)=\hat x_0$ and a.s.~under both $\Q^{\circ,k}$ and $\PP^{\circ,k}$, for every $t\in\R_+$, 
\begin{align}\notag
\hat X^{\circ,k}(t)&=f\left(\Gamma\left[\theta\cdot \left(\hat X^{\circ,k}(0)+\hat W^{\circ,k}(\cdot))+\sigma\hat\Psi^{\circ,k}_A(\cdot)-\sigma_S\hat\Psi^{\circ,k}_S(\cdot)\right)\right](t)\right),\\\notag
\hat Y^{\circ,k}(t)&= \hat X^{\circ,k}(t)- \left(\hat X^{\circ,k}(0)+\hat W^{\circ,k}(\cdot)+\sigma\hat\Psi^{\circ,k}_A(\cdot)-\sigma_S\hat\Psi^{\circ,k}_S(\cdot)\right),
\end{align}

\item $\hat W^{\circ,k}=\hat m+\check A^{\circ,k}-\check D^{\circ,k}$, where 
$(\sigma^{-1}\check A^{\circ,k},\sigma^{-1}_S\check D^{\circ,k})$ is a $2I$-dimensional SBM under $\Q^{\circ,k}$ and
$(\sigma^{-1}\check A^{\circ,k}+\hat\Psi^{\circ,k}_A,\sigma^{-1}_S\check D^{\circ,k}+\hat\Psi^{\circ,k}_A)$ is a $2I$-dimensional SBM under $\PP^{\circ,k}$, both w.r.t.~the filtration 
\begin{align}\notag
\calF^{\circ,k}_t:=\Big\{\check A^{\circ,k}(s),\check D^{\circ,k}(s),\hat W^{\circ,k}(s),\hat X^{\circ,k}(s), \hat Y^{\circ,k}(s), \hat\Psi^{\circ,k}_A(s),\hat\Psi^{\circ,k}_S(s)\;:\; 0\le s\le t\Big\},
\end{align}
\item  
$
 \Psi^{\circ,k}_j(\cdot)=\int_0^\cdot\hat\phi_j(s)ds,$ $j\in\{A,S\},$
for some $[-k,k]^I$-valued, $\calF^{\circ, k}_t$-progressively measurable processes $\{\hat\phi_j=(\hat\phi_{j,i}:i\in[I])\}_j$.
\item 
For every $t\in\R_+,$
\begin{align}\notag
H_{A,i}(t)&=\frac{d(\Q^{\circ,k}\circ (\check A^{\circ,k}_i)^{-1})}{d(\PP^{\circ,k}\circ( \check A^{\circ,k}_i)^{-1})}(t)=\exp\Big\{\int_0^t\hat \phi_{A,i}(s) d\check A^{\circ,k}_{A,i}(s)-\frac{1}{2}\int_0^t \hat \phi^2_{A,i}(s)ds\Big\},
\\\notag
H_{S,i}(t)&=\frac{d(\Q^{\circ,k}\circ (\check D^{\circ,k}_i)^{-1})}{d(\PP^{\circ,k}\circ (\check D^{\circ,k}_i)^{-1})}(t)=\exp\Big\{\int_0^t\hat \phi_{S,i}(s) d\check D^{\circ,k}_{S,i}(s)-\frac{1}{2}\int_0^t \hat \phi^2_{S,i}(s)\rho_ids\Big\}.
\end{align}
\end{enumerate}
\end{lemma}
{\bf Proof.} The $\calC$-tightness and the first three properties follow by standard martingale techniques and Proposition \ref{prop44}. We now prove the forth property. Pay attention that both right-hand sides follow by the second property. Hence, we are only required to establish the left-hand sides, which we provide only for $j=A$.
Fix $t>0$ and a continuous and bounded function $h:\calD[0,T]\to\R$. 
Then, 
\begin{align}\notag
&\E^{\Q^{n,k}}\Big[h\big(\check A^{n,k}_i\big|_{[0,t]}\big)\Big]
=
\E^{\PP^{n,k}}\Big[h\big(\check A^{n,k}_i\big|_{[0,t]}\big)\tfrac{d\Q^{n,k}_{A,i}}{d\PP^{n,k}_{A,i}}\big|_{[0,t]}\Big],
\end{align}
where here and below for any process $E$, its restriction to the time interval $[0,t]$ is denoted by $E|_{[0,t]}$. 
By tightness and converging along a subsequence, we get that
%
%
%
%
$$\E^{\Q^{\circ,k}}\Big[h\big(\check A^{\circ,k}_i\big|_{[0,t]}\big)\Big]
=
\E^{\PP^{\circ,k}}\Big[h\big(\check A^{\circ,k}_i\big|_{[0,t]}\big)H_{A,i}\big|_{[0,t]}\Big].$$ Therefore, $H_{A,i}(t)=d(\Q^{\circ,k}\circ(\check A^{\circ,k}_i)^{-1})/d(\PP^{\circ,k}\circ(\check A^{\circ,k}_i)^{-1}))(t)$. 

\qed

\subsubsection{Convergence of the cost components.}\label{sec434}
Recall Proposition \ref{prop43}. Fix $\eps>0$ and $k_\eps>0$ such that 
\begin{align}\notag
\liminf_{n\to\iy} \;\big\{J(\hat Y^{n,k},\Q^{n,k})-J(\hat Y^{n},\Q^{n})\big\}\ge -\eps.
\end{align}
So, 
\begin{align}\notag
\limsup_{n\to\iy} J(\hat Y^{n},\Q^{n})\le 
\liminf_{n\to\iy}J(\hat Y^{n,k},\Q^{n,k})+\eps\le J(\hat Y^{\circ,k},\Q^{\circ,k})\le V+\eps.
\end{align}
The second inequality follows by the limit given in Lemma \ref{lem45} and the last equality follows by the structure of $\hat Y^{\circ,k}$ and by Proposition \ref{prop31}. This establishes \eqref{4381}.

\vspace{5pt}

\noindent{\bf Acknowledgement.} We thank an anonymous AE and two referees for their suggestions, which helped us to improve our paper.

\footnotesize
\bibliographystyle{abbrv} 
\bibliography{refs} 

\end{document}